\documentclass[12pt]{article}

\usepackage{amssymb}
\usepackage{amsmath}
\usepackage{amsbsy}
\usepackage{amscd}
\usepackage{amsfonts}
\usepackage{amsthm}
\usepackage{mathrsfs}
\usepackage{verbatim}
\usepackage[colorlinks]{hyperref}
\usepackage{fullpage}
\usepackage{mathdots}
\usepackage{graphicx,subfigure}
\usepackage[english]{babel}
\usepackage[utf8]{inputenc}
\usepackage{tikz}
\usepackage{adjustbox}
\usepackage{authblk}
\usepackage{caption}
\usepackage{mathabx}
\usepackage{mathtext}
\usepackage{longtable}
\usepackage{abraces}
\usepackage{stackengine}
\stackMath

\usepackage{tikz-cd}
\usetikzlibrary{backgrounds,fit, matrix}
\usetikzlibrary{positioning}
\usetikzlibrary{calc,through,chains}
\usetikzlibrary{arrows,shapes,snakes,automata, petri}
\usepackage[boxsize=1.25em, centerboxes]{ytableau}

\newtheorem{theorem}{Theorem}[section]
\newtheorem{proposition}[theorem]{Proposition}
\newtheorem{corollary}[theorem]{Corollary}
\newtheorem{lemma}[theorem]{Lemma}

\theoremstyle{definition}
\newtheorem{definition}[theorem]{Definition}
\newtheorem{propdef}[theorem]{Proposition-Definition}
\newtheorem{example}[theorem]{Example}
\newtheorem{notation}[theorem]{Notation}

\theoremstyle{remark}
\newtheorem{remark}[theorem]{Remark}

\numberwithin{equation}{section}

\newcommand{\C}{{\mathbb C}}
\newcommand{\Z}{{\mathbb Z}}

\newcommand{\R}{{\mathbb R}}

\newcommand{\N}{{\mathbb N}}

\newcommand{\vep}{\varepsilon}
\newcommand{\mc}[1]{\mathcal{#1}}
\newcommand{\ms}[1]{\mathscr{#1}}

\DeclareMathOperator{\LL}{\Large{\mathsf{\Lambda}}}

\DeclareMathOperator{\FS}{\stackon[.2pt]{\mathfrak{S}}{\scriptstyle\circ}}
\DeclareMathOperator{\FT}{\stackon[.2pt]{\mathfrak{T}}{\scriptstyle\circ}}
\DeclareMathOperator{\FN}{\stackon[.2pt]{\mathfrak{N}}{\scriptstyle\circ}}

\begin{document}

\title{Locally Semialgebraic Superspaces and Nash Supermanifolds}

\author{Mahir Bilen Can}
\affil{\small{
Tulane University, New Orleans,\\ mahirbilencan@gmail.com}}

\maketitle

\begin{abstract}
In this article, we present a novel theory of locally semialgebraic superspaces along with Nash supermanifolds. By adapting Batchelor's theorem to our framework, we show that all locally semialgebraic superspaces and affine Nash supermanifolds can be derived from appropriate vector bundles. Our analysis of Nash supermanifolds involves an investigation of a topology, where closed sets are defined by Nash subsets. Within this context, we establish the softness of the structure sheaf for an affine Nash manifold. Moreover, we establish that under this topology, the higher cohomology of quasi-coherent sheaves on an affine Nash manifold vanishes completely. These results open up new directions for cohomological studies on Nash manifolds.
\vspace{.5cm}

\\

\end{abstract}

\normalsize

\section{Introduction} 

The main goal of our article is to introduce two new theories of $\Z_2$-graded objects. 
The prototypes of our spaces were developed in the ateliers of mathematical physicists Berezin and Leites~\cite{BerezinLeites}, and the geometer Kostant~\cite{Kostant1975}.
In his seminal work, Kostant expounded, among other things, the prequantizations of $\Z_2$-graded manifolds. 
In particular, he extended the orbit method to the representation theory of Lie supergroups. 
The primary objective of the research was to establish a robust mathematical foundation for the geometry of elementary particles.
In quantum field theories, manifolds must possess the information of both bosonic and fermionic fields, represented by two anti-commuting variables.
Consequently, the formal description of these manifolds necessitates the coordinate functions to assume values in a Grassmann algebra.
Such manifolds are broadly referred to as {\em supermanifolds}~\cite{DeWitt, Rogers1980}. 
The connection between these physical supermanifolds and the $\Z_2$-graded manifolds of Kostant was made precise by Batchelor.
She proved in~\cite{Batchelor1980} that the category of $\Z_2$-graded manifolds is equivalent to the category of (DeWitt) supermanifolds.
From this viewpoint, algebro-geometric versions of supermanifolds make natural appearance in physical theories as well,~\cite{DeligneMorgan1999}.
In fact, this is the underlying theme of Manin's influential book~\cite{Manin} on gauge field theory.
While these exciting mathematical developments took place in the supergravity theory for almost half a century, 
potentially useful semi-algebro-geometric direction has not been pursued. 
Therefore, it is among the goals of our present article to initiate such a study for semialgebraic spaces. 
We achieve this objective by extending the work of Blattner and Rawnsley~\cite{BlattnerRawnsley} to the setting of generalized topologies.
The specific application we have in mind for our theory is a version of an important theorem of Batchelor~\cite{Batchelor1979}, which, roughly speaking, states that every $\Z_2$-graded smooth manifold comes from a vector bundle. 
In fact, we will prove this theorem for two related but distinct categories, that are, 
\begin{enumerate}
\item[(A)] the category of locally semialgebraic spaces, 
\item[(B)] the category of affine Nash manifolds.
\end{enumerate}
We begin our discussion with (A). 
First, we will present a concise overview of the objects of the category of locally semialgebraic spaces, focusing on key concepts that may initially appear complex to the reader. However, we will leave out some details to keep our introduction brief and refer readers to the preliminaries section for more thorough elaboration.
\medskip

Throughout our text, the set of positive integers is denoted by $\Z_+$.
Let $R$ be a real closed field. This means that $R$ is an ordered field with no nontrivial algebraic extension $R_1\supsetneq R$ such that 
$R_1$ is also an ordered field. 
Let $n\in \Z_+$.
A {\em semialgebraic set in $R^n$} is a set that can be presented as a finite union of sets of the form 
\begin{align*}
\{x\in R^n :\ f(x) = 0 , \ g_1(x) > 0, \dots, g_k(x)>0 \},
\end{align*}
where $f,g_1,\dots, g_k$ are some polynomials from $R[x_1,\dots, x_n]$. 
Let $X$ be a semialgebraic set in $R^n$.
A map $f: X\to R$ is called a {\em semialgebraic function} if the graph of $f$ is a semialgebraic subset of $R^{n+1}$. 
Thanks to the ordering of the real closed field $R$, the affine space $R^n$ is equipped with a natural topology for which the open balls form a basis of open subsets. 
We will call this topology the {\it ordinary topology} on $R^n$.
We assume without further notice that all semialgebraic functions we consider in this paper are continuous with respect to the induced ordinary topology on the semialgebraic set under consideration. 
In fact, it is an important property of the semialgebraic functions that they remain continuous with respect to a weaker notion of a topology.
Since this fact is an essential ingredient in our results, and since it further motivates our discussion, we will explain our remark here.
The {\em restricted topology}, in the sense of ~\cite[\S 7]{DelfsKnebusch1981II}, is defined as follows.
Let $Y$ be a set. 
A family $\FS(Y)$ of subsets of $Y$ is called a {\em restricted topology on $Y$} 
if the following two conditions are satisfied: 
\begin{enumerate}
\item $\{\emptyset,Y\}\subseteq \FS(Y)$, 
\item $\FS(Y)$ is closed under finite intersections and finite unions of its elements. 
\end{enumerate}
Then a function $f:X\to Y$ between two restricted topological spaces is called {\em continuous} if for every $U\in \FS(Y)$,
the preimage $f^{-1}(U)$ is an element of $\FS(X)$. 
\medskip

In a semialgebraic set $X$, the open semialgebraic subsets form a restricted topology.
In this case, the pair $(X,\ms{O}_X)$, where $\ms{O}_X$ is the sheaf of semialgebraic functions on $X$ (with respect to the restricted topology)
becomes a ringed space over $R$. 
An {\em affine semialgebraic space} is a ringed space over $R$ that is isomorphic to the ringed space over $R$ of a semialgebraic 
set in $R^n$, for some $n\in \Z_+$. 
Among the examples of affine semialgebraic spaces are the real algebraic varieties (projective or not).
More generally, a {\em semialgebraic space} is a ringed space that has a finite open cover by affine semialgebraic subspaces. 
\medskip

Even though a restricted topology is not a true topology, it has proven to be valuable in resolving issues within semialgebraic geometry. Delfs and Knebusch have introduced in~\cite{DelfsKnebusch} a stricter concept, known as a {\em generalized topological space}, which leads to a substantial expansion of the theory of semialgebraic spaces.
The precise definition of a generalized topology is more technical than the definition of a restricted topology. 
We postpone it to the preliminaries section.
A {\em locally semialgebraic space} is a ringed space $(X,\ms{O}_X)$ over $R$,
where $X$ is a generalized topological space that is equipped with an admissible covering by open semialgebraic subspaces.
Hereafter, if there is no danger for confusion, we denote a locally semialgebraic space $(X,\mathscr{O}_X)$ by $X$ only. 
This paper will concentrate on locally semialgebraic spaces that are both regular and paracompact.
To keep things brief in the introduction, we defer the precise definitions of these technical notions to the preliminaries section. 
It is important to remember that every semialgebraic space is inherently a locally semialgebraic space. Additionally, it should be noted that any semialgebraic subspace of a regular locally semialgebraic space is affine.

\medskip

Let us write $\overline{0}$ (resp. $\overline{1}$) for the identity (resp. non-identity) element of $\Z_2$. 
By a {\em $\Z_2$-graded ring}, we mean an associative ring with unity, denoted $A$, together with a direct-sum decomposition of its additive group structure, $A= A_{\overline{0}} \oplus A_{\overline{1}}$, such that $A_{m} A_{n} \subset A_{m+n}$ for $\{m,n\}\subset \Z_2$. 
In this case, for $m\in \{\overline{0},\overline{1}\}$, an element $a$ from $A_m$ is said to be a {\em homogeneous element of degree $m$.}
As before, let $R$ denote a real closed field. 
Let us assume that $A$ has an $R$-algebra structure as well.
If for every $a\in A_m$ and $b\in A_n$, where $\{m,n\}\subset \{\overline{0},\overline{1}\}$, the equality $ab = (-1)^{m+n} ba$ holds, then $A$ is said to be a {\em $\Z_2$-graded commutative $R$-algebra}. 
Let $X$ be a locally semialgebraic space.
Let $\ms{A}$ be a sheaf of $\Z_2$-graded commutative $R$-algebras on $X$.
Let $s\in \Z_+$.
We call the pair $(X,\ms{A})$ a {\em locally semialgebraic superspace of odd-dimension $s$} if the following two conditions hold:
\begin{enumerate}
\item There is a short exact sequence of sheaves of $\Z_2$-graded commutative algebras of the form,
\begin{align*}
0\to \ms{A}' \to \ms{A} \to \mathscr{O}_X \to 0,
\end{align*}
where the $\Z_2$-grading on $\mathscr{O}_X$ is given by $(\mathscr{O}_X)_{\overline{0}}:=\mathscr{O}_X$ and $(\mathscr{O}_X)_{\overline{1}}:=0$.
\item There is an admissible cover by open semialgebraic spaces $\bigcup_{i\in I} U_i = X$,   
where, for each $i\in I$, the restricted sheaf $\ms{A}|_{U_i}$ is isomorphic to the restricted sheaf of $\Z_2$-graded commutative algebras, 
$(\mathscr{O}_X \otimes \LL R^s)|_{U_i}$. 
Here, $\mathscr{O}_X$ is $\Z_2$-graded as in part 1, and $\LL R^s$ is the total exterior algebra of $R^s$, whose $\Z_2$-grading  
is determined by the subspaces
\begin{align}\label{A:decompositionoftotal}
(\LL  R^s)_{\overline{0}} := \bigoplus_{i\ : \ \text{even}} \LL^i  R^s \qquad \text{and}\qquad (\LL  R^s)_{\overline{1}} := \bigoplus_{i\ : \ \text{odd}} \LL^i R^s
\end{align}
\end{enumerate}

It is not difficult to adapt the usual definition of a topological vector bundle to the setting of the locally semialgebraic spaces. 
In fact, later in this text we will work with even more general notion of a vector bundle that is defined on a generalized topological space. 
For now, we proceed with the assumption that $p: Y\to X$ is a locally semialgebraic vector bundle. 
In this case, the sum of all exterior powers of $Y$, denoted $\LL Y$, is a locally semialgebraic vector bundle also. 
Moreover, it is easily seen that the space of global sections of $\LL Y$, denoted $\Gamma(X,\LL Y)$, has the structure of a 
$\Z_2$-graded commutative algebra, where the $\Z_2$-grading is defined by 
\begin{align*}
\Gamma(X,\LL Y)_{\overline{0}} := \Gamma\left(Y\ , \bigoplus_{i\ : \ \text{even}} \LL^i Y\right) \qquad \text{and}\qquad 
\Gamma(X,\LL Y)_{\overline{1}} := \Gamma\left(Y\ , \bigoplus_{i\ : \ \text{odd}} \LL^i Y\right).
\end{align*}
In this notation, the thrust of our paper is the following result, which is analogous to the main theorem of~\cite{Batchelor1979}.

\begin{theorem}\label{T1:intro}
Let $X$ be a regular, paracompact locally semialgebraic space. 
If $(X,\ms{A})$ is a locally semialgebraic superspace of odd-dimension $s$, 
then there exists a locally semialgebraic vector bundle $E$ such that $\ms{A}(X)$ is isomorphic to $\Gamma(X,\LL E)$ as a $\Z_2$-graded commutative algebra. 
\end{theorem}

The proof of Theorem~\ref{T1:intro} is interesting by itself. 
Indeed, as we mentioned before, we introduce a metatheory. 
This theory yields immediately a description of the potential obstruction to the existence of a vector bundle $E$ as in our Theorem~\ref{T1:intro}. 
Then we show that, under our assumptions, there is no obstruction, namely the vanishing of higher cohomology, for certain sheaves of $\ms{O}_X$-modules that are related to our vector bundle.
\medskip

We now proceed to introduce the ingredients of our second main result. 
To this end, we specialize our real closed field $R$ to the field of real numbers, $\R$.  
Let $n\in \Z_+$.
Roughly speaking, an {\em affine Nash submanifold in $\R^n$} is a semialgebraic, analytic submanifold of $\R^n$.
In classical algebraic geometry, an abstract variety is defined by gluing affine algebraic varieties. 
Similarly to the construction of an abstract algebraic variety, an abstract {\em Nash manifold} is obtained by gluing affine Nash manifolds. 
There are deep connections between Nash manifolds and real algebraic varieties. 
To give a general example let us consider an orbit of a real algebraic group action.
The orbit is not necessarily a real algebraic variety but it always has the structure of a Nash manifold. 
We will give a more specific example of this phenomenon later in this text. 
It should be noted that when considering algebraic group actions defined over algebraically closed fields, there are no issues arising from objects outside the relevant algebraic category. In other words, every orbit resulting from an algebraic group action has a quasi-projective variety structure that is intrinsic to it.
\medskip

Let $M$ be a Nash manifold. 
There are several useful topologies on $M$.  
The first one is the topology that is inherited from the euclidean topology of $\R^n$. 
The second one is the restricted topology whose open sets are given by the open semialgebraic subsets of $M$. 
Now, for our second main result, we work with a weaker topology, which is similar to the Zariski topology on the prime spectrum of a ring. 
Let $f: U\to \R$ be a semialgebraic function defined on an open semialgebraic subset $U\subseteq M$. 
If, in addition, $f$ is a real analytic function, then $f$ is called a {\em Nash function} on $U$.  
It turns out that, under the usual point-wise addition and multiplication of functions, and the scaling action of $\R$ on them, 
the Nash functions on $U$ form a noetherian $\R$-algebra.
Our notation for this algebra is $\ms{N}_M(U)$. 
\medskip

We call a subset $V\subseteq M$ a {\em Nash subset} if $V$ is given by the common zero-set of a finite collection of global Nash functions on $M$. 
The set of all Nash subsets of $M$ is closed under arbitrary intersections and finite unions. 
This fact follows from~\cite[Proposition 8.6.2]{BochnakCosteRoy}.
Thus, the collection of all Nash subsets of $M$ defines a topology, which we call here the {\em $\mathfrak{N}$-topology of $M$.}
When we want to emphasize this topology, we will write $(M,\mathfrak{N})$ instead of $M$ only.
The subtlety that we want to point out with this notion is that, in our next theorem, we work with the sheaves that are defined on the $\mathfrak{N}$-topology of $M$. 
We call the open subsets (resp. the closed subsets) of an $\mathfrak{N}$-topology the {\em $\mathfrak{N}$-open subsets} (resp. the $\mathfrak{N}$-closed subsets). 
While moving forward, our permanent assumption in this text will be that all of our $\mathfrak{N}$-topologies are {\em separated}.
This means that the diagonal copy of $M$ in $M\times M$ is an $\mathfrak{N}$-closed subset of $M\times M$.
Note that the property of being separated does not necessarily imply the Hausdorff property, as the $\mathfrak{N}$-topology on $M \times M$ may not always coincide with the product topology.
\medskip

Let $(M,\mathfrak{N})$ be an $n$-dimensional Nash manifold.  
We view $\ms{N}_M$ as a presheaf on $(M,\mathfrak{N})$.
In fact, it is not difficult to check that $\ms{N}_M$ is a sheaf of $R$-algebras on $(M,\mathfrak{N})$.
Let $s\in \Z_+$.
We call a pair $(M,\ms{A})$, where $\ms{A}$ is a sheaf of $\Z_2$-graded commutative $\R$-algebras on $(M,\mathfrak{N})$, 
a {\em Nash supermanifold of odd-dimension $s$} if the following conditions hold:
\begin{enumerate}
\item There is a short exact sequence of sheaves $\Z_2$-graded commutative algebras of the form,
\begin{align*}
0\to \ms{A}' \to \ms{A} \to \ms{N}_M \to 0,
\end{align*}
where, the $\Z_2$-grading on $\ms{N}_M$ is defined by setting 
\begin{align*}
(\ms{N}_M)_{\overline{0}}:=\ms{N}_M\quad \text{ and }\quad (\ms{N}_M)_{\overline{1}}:=0.
\end{align*}

\item There is a finite open cover (with respect to the $\mathfrak{N}$-topology), $\bigcup_{i\in I} U_i = M$, 
where, for every $i\in I$, the restricted sheaf $\ms{A}|_{U_i}$ on $U_i$ is isomorphic to the restricted sheaf of $\Z_2$-graded commutative algebras, 
$(\ms{N}_M \otimes \LL \R^s)|_{U_i}$. 
Here, $\ms{N}_M$ is $\Z_2$-graded as in part 1, and $\LL \R^s$ is the total exterior algebra of the euclidean space $\R^s$, whose grading  
is determined by the subspaces,
\begin{align*}
(\LL \R^s)_{\overline{0}} := \bigoplus_{i\ : \ \text{even}} \LL^i \R^s \qquad \text{and}\qquad (\LL \R^s)_{\overline{1}} := \bigoplus_{i\ : \ \text{odd}} \LL^i \R^s
\end{align*}
\end{enumerate}
Let $\R[[ X_1,\dots, X_n ]]_{\text{alg}}$ denote the trivially graded ring of power series which are algebraic over the polynomial ring $\R[X_1,\dots, X_n]$. 
In a nutshell, the above definitions are saying that a Nash supermanifold of odd-dimension $s$ 
is a Nash manifold together with a sheaf of $\Z_2$-graded commutative algebras 
that is locally isomorphic to a sheaf whose stalk at $x\in M$ is isomorphic to the $\Z_2$-graded commutative algebra,
\begin{align*}
\LL \R^s \otimes_\R \R[[ X_1,\dots, X_n ]]_{\text{alg}}.
\end{align*}
We proceed to fulfill the natural desire of determining the global structure of an affine Nash supermanifold.

As before, let $(M,\mathfrak{N})$ be an (affine) Nash manifold. 
A vector bundle $p: Y\to M$ is called an {\em (affine) Nash vector bundle} on $M$ if $Y$ is an (affine) Nash manifold such that 1) $p$ is a Nash mapping, 2) for every $x\in M$, the addition and the scalar multiplication on the vector space $p^{-1}(x)$ are Nash mappings. 
Now let $Y$ be an (affine) Nash vector bundle on $M$. 
Then the sum of all exterior powers of $Y$, denoted $\LL Y$, is an (affine) Nash vector bundle on $M$ as well. 
The space of global sections of $\LL Y$ on $(M,\mathfrak{N})$ has the structure of a $\Z_2$-graded commutative algebra.
In this notation, our second main result is the following statement. 

\begin{theorem}\label{T2:intro}
Let $(M,\ms{A})$ be an affine Nash supermanifold of odd-dimension $s$. 
Then there exists an affine Nash vector bundle $E$ on $M$ such that $\ms{A}(M)$ is isomorphic, 
as a $\Z_2$-graded commutative algebra, to $\Gamma(M,\LL E)$. 
\end{theorem}

The main idea behind the proof of Theorem~\ref{T2:intro} is similar to the proof of Theorem~\ref{T1:intro}.
Nevertheless, for completing the proof, we utilized some additional facts about the sheaf of Nash functions.
A sheaf $\ms{F}$ on a topological space $X$ is called {\em flabby} (respectively, {\em soft}) 
if for every open (respectively, closed) subset $V\subseteq X$, the restriction map $\Gamma (X,\ms{F})\to \Gamma(V,\ms{F})$, $s\mapsto s\vert_V$ is surjective. 
A flabby sheaf, also called a {\em flasque sheaf}, on a paracompact space is a soft sheaf,~\cite[Theorem 3.3.1]{Godement}.
The proof of the next standard fact that we want to point out can be deduced from~\cite[Theorem 4.4.3]{Godement}: 
A soft sheaf on a paracompact topological space $X$ is acyclic if the underlying topology of $X$ is Hausdorff.
It is important here to recall that while the $\mathfrak{N}$-topology on a Nash manifold is separated according to our convention, it is not Hausdorff unless $M$ is a finite set. 
We should also mention here that every $\mathfrak{N}$-topology is paracompact.
For the purpose of our discussion, we record the following assertion regarding the softness of $\mathcal{N}_M$.

\begin{theorem}\label{T:Nashsoft}
Let $M$ be an affine Nash manifold. Then the sheaf $\ms{N}_M$ is a soft sheaf with respect to the $\mathfrak{N}$-topology on $M$.
\end{theorem}
The proof of Theorem~\ref{T:Nashsoft} is obtained through a straightforward application of the {\em Efroymson's extension theorem}, which is proven in the articles~\cite{Efroymson, Pecker}.
Although our sheaf of Nash functions $\ms{N}_M$ is soft, it fails to be a flabby sheaf. 
Fortunately, the category of sheaves of $\ms{N}_M$-modules that we are interested in exhibits the desired outcome of cohomology vanishing, which flabbiness could provide.

Let $M$ be an affine Nash manifold. 
Let $A$ denote the $\R$-algebra of all Nash functions on $M$. 
For $f\in A$, let $E(f)$ denote the $\mathfrak{N}$-open subset defined by $$E(f):=\{x\in M:\ f(x)\neq 0\}.$$
A sheaf of $\ms{N}_M$-modules $\ms{F}$ is said to be a {\em quasi-coherent sheaf of $\ms{N}_M$-modules} 
if there exists an $A$-module $K$ such that, for every $f\in A$, $\ms{F}(E(f))$ is the localization of the ring $A$ at the element $f$.
We will discuss these sheaves of modules in detail in Section~\ref{S:Affinenash}.
For now let us mention that the quasi-coherent sheaves of $\ms{N}_M$-modules on $M$ give the quasi-coherent sheaves of $\ms{O}_X$-modules on the affine scheme $X:=\text{Spec}(A)$. 
The main technical result that we will use for proving Theorem~\ref{T2:intro} is the following vanishing result, which is well-known in the context of affine schemes. 
\begin{theorem}\label{T3:intro}
Let $M$ be an affine Nash manifold. Let $\ms{F}$ be a quasi-coherent sheaf of $\ms{N}_M$-modules on $M$.
Then we have $H^i(M,\ms{F}) = 0$ for all $i>0$.
\end{theorem}

We are now ready to describe the structure of our paper. 
The foundations of locally semialgebraic spaces, in the spirit of modern algebraic geometry, were laid by Delfs and Knebusch
in several installments. 
Initially, they formulated the theory of semialgebraic spaces by embedding them within the real points of real algebraic varieties. Subsequently, they produced a significant manuscript that focused on locally semialgebraic spaces.
After that Delfs continued to develop homological methods for the locally semialgebraic spaces. 
In the preliminaries section (Section~\ref{S:Preliminaries}) we make a brief introduction to the work of Delfs and Knebusch.
Among the important results that we will review there are the {\em theorem of partition of unity} and the {\em extension theorem} of 
locally semialgebraic functions. These results will play a crucial role in proving our theorems.  
In~\cite{BlattnerRawnsley}, Blattner and Rawnsley describe a far reaching generalization of Bathchelor's theorem. 
The gist of their work is to reduce Batchelor's theorem to a vanishing statement for the first cohomology groups of certain sheaves. 
After we setup some terminology, in Section~\ref{S:Ringed}, we observe that the methods of Blattner and Rawnsley can be adapted in the general setting of a generalized topology.
Although it looks somewhat abstract, the generality we achieve by this adaptation gives us enough space to prove Batchelor's Theorem 
for the regular locally semialgebraic spaces in Section~\ref{S:Locally}.
We give a proof of Theorem~\ref{T2:intro} in Section~\ref{S:Affinenash}.
As expected, the proof outline of this theorem shares similarities with that of Theorem~\ref{T1:intro}. An integral element of the proof, namely Theorem~\ref{T3:intro} concerning the vanishing of higher cohomology sheaves, is also proven in Section~\ref{S:Affinenash}.
We close our paper by Section~\ref{S:Final}, where we mention our future work about the group objects in the category of Nash supermanifolds.

\section{Preliminaries}\label{S:Preliminaries}

In this article, we break new ground for certain noncommutative semialgebraic geometry.
The key notion for our development is the notion of a semialgebraic space.

Not only in algebraic geometry but also in semialgebraic geometry, the finiteness of the coverings is a built-in property. 
In this sense, the notion of a restricted topology, which we introduced earlier, is a natural choice for our purposes. 
Recall that a family $\FS(X)$ of subsets of a set $X$ is called a restricted topology on $X$ if it satisfies the following conditions:
\begin{itemize}
\item $\{\emptyset, X \} \subseteq \FS(X)$, 
\item $\FS(X)$ is closed under finite intersections and finite unions.
\end{itemize}
In this case $(X,\FS(X))$ is called a {\em restricted topological space}. 
The elements of $\FS(X)$ are called the {\em open sets} of the restricted topology. 
A {\em continuous map} of restricted topological spaces is a function, $f: X\to Y$, such that for every $U\in \FS(Y)$,
the set $f^{-1}(U)$ is open in $X$.
From now on, for the sake of simplicity in notation, we will adopt our earlier convention: the restricted topological space $(X,\FS(X))$ will be referred to simply as $X$.
\medskip

Let us now introduce a fundamental example of a restricted topological space for our intended purposes.
Let $n\in \Z_+$. Let $R$ be a real closed field. 
The collection of all semialgebraic subsets of $R^n$ is the smallest family of subsets containing all sets of the form 
\begin{align}\label{A:basicopen}
\{x\in R^n :\ f(x) > 0\},\qquad \text{where $f\in R[x_1,\dots, x_n]$},
\end{align}
and closed under finite unions, finite intersections, and complements. 
We call a semialgebraic subset of $R^n$ an {\em open semialgebraic subset} if it is open in the ordinary topology on $R^n$. 
A {\em closed semialgebraic subset of $R^n$} is the complement of an open semialgebraic subset in $R^n$.
\medskip

The dimension of a semialgebraic subset $X\subseteq R^n$ can be defined in purely algebraic terms as follows.
Let $\mc{I}(X)$ denote the ideal of polynomials $f\in R[x_1,\dots, x_n]$ such that $f(x) = 0$ for every $x\in X$. 
Then the {\em dimension of $X$} is defined to be the Krull dimension of $R[x_1,\dots, x_n]/\mc{I}(X)$.  
Note that the dimension of a semialgebraic subset $X\subseteq R^n$ is given by the dimension of the Zariski closure of $X$ in $R^n$.

\begin{example}\label{E:threeexamples}
The semialgebraic subsets of $\R^1$ are exactly the finite unions of open intervals and points. 
Let us consider the following subsets of $\R^1$:
\begin{itemize}
\item $X_1:= \{ x\in \R^1 :\ x \geq 0\}$,
\item $X_2:=\{ x\in \R^1 :\ 0< x <1\} \cup \{ x\in \R^1 :\ 1< x <\infty\}$,
\item $X_3:=\Z$.
\end{itemize}
Then $X_1$ is a one dimensional closed semialgebraic subset of $\R^1$, and $X_2$ is a one dimensional open semialgebraic subset of $\R^1$.
However, $X_3$ is not a semialgebraic subset of $\R^1$. 
Notice also that the Zariski closures of all of these subsets in $\R^1$ are equal to $\R^1$.
\end{example}

It is easy to check that the collection of all open semialgebraic subsets of $\R^n$, including the empty set, is a restricted topology on $\R^n$. 
Let $X\subseteq \R^n$. 
By intersecting the open semialgebraic subsets of $\R^n$ with $X$, we obtain the open sets of a restricted topology on $X$.
We will refer to this restricted topology the {\em $\mathfrak{S}$-topology on $X$}.
\medskip

Restricted topologies, as well as the upcoming discussion on generalized topologies, illustrate a more abstract concept of a Grothendieck site.  
However, in this paper, we will not extensively explore the theory of sheaves on Grothendieck sites as presented in~\cite{Artin} or~\cite{Knutson}. 
Nonetheless, we find a relevant resource in \cite[Section 3.1]{EdmundoPrelli}, 
which provides a comprehensive overview of the current state of art concerning sheaves on the Grothendieck site of a locally semialgebraic space.
\medskip

All presheaf operations on the categories of (pre)sheaves on restricted topological spaces are defined analogously 
to the presheaf operations on the categories of (pre)sheaves on the ordinary topological spaces. 
For example, if $\mathscr{F}$ and $\mathscr{G}$ are two presheaves on a restricted topology on $X$, 
then a presheaf map $\mathscr{F}\to \mathscr{G}$ is defined 
by the requirement that for every inclusion $U\hookrightarrow V$ in $\FS(X)$, there is a commutative diagram,
\[
\begin{tikzcd}
  \mathscr{F}(V) \arrow[d] \arrow[r] &  \mathscr{G}(V) \arrow[r, d] \\
 \mathscr{F}(U) \arrow [r]  &  \mathscr{G}(U)
\end{tikzcd}
\]

Shortly, we will engage with sheaves of $R$-algebras on locally semialgebraic spaces, where $R$ is an appropriately defined field. Such a space comes with a generalized topology (to be defined) on it. 
In preparation for explaining these notion, in our next definition, $X$ will represent either a restricted topological space, a standard topological space, or a generalized topological space.

\begin{definition}
A {\em ringed space over $R$ on $X$} is a pair of the form $(X,\ms{A})$, where $\ms{A}$ is a sheaf of $R$-algebras on $X$.
If every stalk of $\ms{A}$ is a local ring, then we call $(X,\ms{A})$ a {\em locally ringed space over $R$ on $X$}. 

Let $(X,\ms{A})$ and $(Y,\ms{B})$ be two ringed spaces over $R$.
A morphism between them is a pair $(\varphi , \vartheta)$, where $\varphi : X\to Y$ is a continuous map, 
and $\vartheta : \ms{B} \to \varphi_* \ms{A}$ is a sheaf map.
Here, $ \varphi_* \ms{A}$ is the {\em direct image sheaf} on $Y$, 
defined by $ (\varphi_* \ms{A}) (U) := \ms{A}(\varphi^{-1}(U))$ where $U$ is an object (that is, an open set) in $Y$. 
\end{definition}

The general theory of semialgebraic spaces is built on real closed fields. 
Recall that a {\em real closed field} is an ordered field that is maximal with respect to inclusions of ordered fields. 
Equivalently, by~\cite[Theorem 1.2.2]{BochnakCosteRoy}, a real closed field is a non-algebraically closed field $R$ such that 
the degree two extension of $R$ that is obtained by adding the solution of the equation $x^2 +1 = 0$ to it is an algebraically closed field. 
We proceed with the assumption that $R$ is a real closed field.

Recall also that an affine semialgebraic space is a ringed space over $R$ that is isomorphic to $(X,\mathscr{O}_X)$, 
where $X$ is a semialgebraic set in $R^n$ for some $n\in \Z_+$, and $\mathscr{O}_X$ is the sheaf of semialgebraic functions on $X$. 

A {\em semialgebraic subset} of an affine semialgebraic space $(X,\ms{O}_X)$ is a subset of $X$ 
For such an affine semialgebraic space, a subset $A$ is said to be {\em semialgebraic} 
if the image of $A$ in $R^n$ is isomorphic to a semialgebraic set in $R^n$ also.  
A {\em semialgebraic space} is a ringed space $(X,\mathscr{O}_X)$ over $R$ that has a finite open cover
$X=\bigcup_{i\in I} X_i$ such that each ringed space $(X_i , \mathscr{O}_X \vert_{X_i})$, where $i\in I$,  is an affine semialgebraic space.
In particular, the sheaf $\mathscr{O}_X$ can be viewed as a subsheaf of the sheaf of all $R$-valued functions on $X$. 
Now, let $(X,\mathscr{O}_X)$ be a semialgebraic space over $R$.
Let $X=\bigcup_{i\in I} X_i$ be a finite open cover of $X$, where the $(X_i , \mathscr{O}_X \vert_{X_i})$'s are affine semialgebraic spaces. 
Then, for each $i\in I$, the ringed space $(X_i , \mathscr{O}_X \vert_{X_i})$ is isomorphic to a semialgebraic set in $R^n$ for some $n\in \Z_+$. 
A subset $A\subseteq X$ is called a {\em semialgebraic subset} if for every $i\in I$, the intersection $A\cap X_i$ is a semialgebraic subset of the affine semialgebraic space $(X_i,\ms{O}_X |_{X_i})$.

In his work \cite{Delfs}, Delfs developed a sheaf cohomology theory for semialgebraic geometry. 
Before delving into Delfs' theorems, we will first provide a comprehensive introduction to the concept of a locally semialgebraic space.

Let $X$ be a set. A {\em generalized topology} on $X$ is a set $\FT(X)$ of subsets of $X$, 
called the {\em open subsets of $X$}, and a set $\text{Cov}_X$ of families $\{ U_i  | i\in I \}$ in $\FT(X)$, called {\em admissible coverings},
such that the following properties hold:
\begin{enumerate}
\item[1.] The empty set and $X$ itself is an element of $\FT(X)$. 
\item[2.] All finite unions and intersections of the elements of $\FT(X)$ are also elements of $\FT(X)$. 
\item[3.] Every family $\{ U_i  |   i\in I \}$, where $I$ is a finite set, is an element of $\text{Cov}_X$.
\item[4.] For every $\{ U_i | i\in I \} \in \text{Cov}_X$, the union $\bigcup_{i\in I} U_i$ is an element of $\FT(X)$. 
\end{enumerate} 
For $U\in \FT(X)$, the set of all coverings $\{U_i | i\in I\}$ such that $\bigcup_{i\in I} U_i = U$ is denoted by $\text{Cov}_X(U)$.
The elements of $\text{Cov}_X(U)$ are called the {\em admissible coverings of $U$}.
\begin{enumerate}
\item[5.] Let $U\in \FT(X)$. 
If $\{ U_i  | i\in I \} \in \text{Cov}_X(U)$ and $V\in \FT(X)$, then we have 
\[
\{ V\cap U_i | i\in I \} \in \text{Cov}_X(V).
\]
\item[6.] Let $U\in \FT(X)$. 
If $\{ U_i | i\in I \} \in \text{Cov}_X(U)$ and $\{ V_{ij}  | j \in I_i \} \in \text{Cov}_X(U_i)$, then we have 
\[
\{ V_{ij} | i\in I, j\in I_i \} \in \text{Cov}_X(U).
\]
\item[7.] If $\{ U_i | i\in I \}$ is a family in $\FT(X)$ such that $U:=\bigcup_{i\in I} U_i$ is an element of $\FT(X)$,
and if $\{V_j | j\in J\}\in \text{Cox}_X(U)$ is a refinement of $\{ U_i | i\in I \}$, that is to say, for every $j\in J$, there exists $i\in I$ such that $V_j \subseteq U_i$, 
then we have 
\[
\{ U_i | i\in I \} \in \text{Cov}_X(U).
\] 
\item[8.] Let $U\in \FT(X)$. 
Let $\{ U_i | i\in I \} \in \text{Cov}_X(U)$. If $V$ is a subset of $U$ such that $U_i \cap V \in \FT(X)$ for every $i\in I$, then $V\in \FT(X)$.
In particular, by 5., we have $\{ U_i \cap V | i\in I \} \in \text{Cov}_X(V)$.
\end{enumerate}

\begin{definition}\label{D:lsa}
A {\em locally semialgebraic space} is a ringed space $(X,\ms{O}_X)$ over $R$, where $X$ is a generalized topological space 
having an admissible covering $\{U_i   | i\in I \} \in \text{Cov}_X$ such that the ringed subspace $(U_i, \ms{O}_X \vert_{U_i})$ over $R$ 
is a semialgebraic space for every $i\in I$.
A morphism of locally semialgebraic spaces is a morphism of the underlying ringed spaces. 
\end{definition}

\begin{remark}
\begin{enumerate} 
\item 
Clearly, if $(X,\ms{O}_X)$ is an ordinary semialgebraic space, then its restricted topology $\FS(X)$ satisfies all the axioms of a generalized topology. 
In other words, every semialgebraic space is automatically a locally semialgebraic space. 
More generally, let us consider a restricted topological space, $(X,\ms{O}_X)$.
Let us define $\FT(X) := \FS(X)$. 
For the coverings of an element $U\in \FT(X)$, we consider the set of all families $\{U_i | i\in I\}$ in $\FS(X)$ such that $\bigcup_{i\in I} U_i = U$,
and $U$ is already covered by finitely many of the $U_i$'s ($i\in I$). 
Then it is easy to check that with these coverings and open sets, $(X,\ms{O}_X)$ is equipped with a generalized topological space structure. 
In other words, every restricted topological space can be regarded as a generalized topological space. 
Thus, any general theorem that is applicable to the generalized topological spaces is applicable to the  
restricted topological spaces as well. 

\item The ringed space of a locally semialgebraic space is actually a {\em locally ringed space}, that is, for every $x\in X$, the stalk $\ms{O}_{X,x}$ is a local ring. 
In fact, its quotient ring $\ms{O}_{X,x}/\mathfrak{m}_x$, where $\mathfrak{m}_x$ is the unique maximal ideal of $\ms{O}_{X,x}$, is isomorphic to the real closed field, $R$. 
\end{enumerate}
\end{remark}

A concise but very useful summary of the theory of locally semialgebraic spaces, up to year 1984, can be found in the article~\cite{DelfsKnebusch1984} of Delfs and Knebusch.
\medskip

For our next definition, we fix a locally semialgebraic space $M$ over $R$. 
This means that we fix an admissible covering $\{M_i | i \in I\} \in \text{Cov}_M$, where each ringed subspace $(M_i,\ms{O}_M |_{M_i})$, $i\in I$, is a semialgebraic space.

\begin{definition}
We maintain our notation from the previous paragraph. 
A {\em locally semialgebraic subset} of $M$ is a subset $X\subseteq M$ such that for every element $W\in \FT(M)$ where $(W,\ms{O}_M|_W)$ is a semialgebraic space, the intersection $X\cap W$ is a semialgebraic subset of the semialgebraic space $W$. 
The set of all locally semialgebraic subsets of $M$ is denoted by $\mathfrak{T}(M)$.
\end{definition}

\begin{remark}\label{R:tocheckthat}
To check that a subset $X\subseteq M$ is a locally semialgebraic subset, it suffices to check that $X\cap M_i$ is a semialgebraic subset of $M_i$ for every $i\in I$. 
\end{remark}

\begin{example}\label{E:U_r's}
Let $\FT(\R^2)$ denote the set $\{ U_r :\ r\in \N\cup \{\infty\} \}$, where $U_r$ is defined by
\[
U_r :=
\begin{cases}
\emptyset & \text{ if $r=0$},\\
\{x\in \R^2 :\ |x| < r \} & \text{ if $r\in \Z_+$},\\
\R^2 & \text{ if $r=\infty$}.
\end{cases}
\]
It is easy to check that the set of all families of the form $\{ U_r :\ r\in I\}$, where $I\subseteq \N\cup \{\infty\}$, satisfies all axioms of being a set of admissible coverings. Let us denote this set of families by $\text{Cov}_{\R^2}$.
Hence, the pair $(\FT(\R^2),\text{Cov}_{\R^2})$ gives a generalized topology on $\R^2$. 
For every $r\in \Z_+$, we will view $U_r$ as an open semialgebraic subspace of the standard affine semialgebraic space $(\R^2,\ms{O}_{\R^2})$. 
Then we define $(X,\ms{O}_X)$ as the inductive limit $X:=\varinjlim_{r\in \N} U_r$ in the category of locally semialgebraic spaces. 
The fact that this inductive limit exists in the category of locally semialgebraic spaces follows from~\cite[Lemma 1.2]{DelfsKnebusch1984}.
Notice that an admissible covering of $X$ is given by $\{ U_r :\ r\in \N\} \in \text{Cov}_{\R^2}$.

While on a set-theoretical level, $X$ can be readily identified with $\mathbb{R}^2$, as a locally semialgebraic space it is different from the standard affine semialgebraic space $(\mathbb{R}^2, \ms{O}_{\mathbb{R}^2})$. In fact, $(X, \ms{O}_X)$ does not even qualify as a semialgebraic space. This becomes evident when we consider the set of closed points in the generalized topology of $(X, \ms{O}_X)$, which is found to be empty. This follows from the fact that the basis of the generalized topology of $X$ consists of the open sets $U_r$, where $r \in \mathbb{N} \cup {\infty}$. A closed semialgebraic subset $Y$ of $X$ can either be empty or unbounded. However, it is easy to see that every nonempty semialgebraic set possesses a closed point. Consequently, $(X, \ms{O}_X)$ is not a semialgebraic space.
\end{example}
\medskip

Two other important properties of the (set of) locally semialgebraic subsets are as follows. 
The proofs of these properties are not difficult. 
\begin{enumerate}
\item If $X$ is an element of $\mathfrak{T}(M)$, then so is $M\setminus X$.
\item Both the union and the intersection of any locally finite family from $\mathfrak{T}(M)$ are contained in $\mathfrak{T}(M)$. 
\end{enumerate}

Let us explain what we mean by a locally finite family in this context. 
\begin{definition}\label{D:locallyfinite}
Let $M$ be a generalized topological space. 
A family $\{X_j | j\in J\}$ of subsets of $M$ is called {\em locally finite} if any element $W\in \FT(M)$ meets only finitely many $X_j$'s.
Equivalently, for every $i\in I$, the open set $M_i$ meets only finitely many $X_j$'s. 
We note in passing that every locally finite family in $\FT(M)$ is an element of $\text{Cox}_M$.
\end{definition}

\medskip

So far, we dealt with locally semialgebraic subsets. 
Next, we want to explain the notion of a locally semialgebraic subspace.
\medskip

Let $M$ and $\mathfrak{T}(M)$ denote, as before, a locally semialgebraic space and the set of locally semialgebraic subsets of $M$, respectively. 
Let $\{ M_i | i\in I\}$ be an admissible covering of $M$. 
We enlarge this admissible covering by adding to it all finite unions of its elements. 
Let $\{ N_j | j\in J\}$ denote the resulting family of open semialgebraic subsets of $M$. 
It is easily checked that $\{ N_j | j\in J\}$ is an admissible covering of $M$ as well. 
Let $X$ be an element of $\mathfrak{T}(M)$. 
Then every intersection $N_j \cap X$ ($j\in J$) is a semialgebraic subspace of $N_j$. 
\begin{propdef}
We maintain our notation from the previous paragraph. 
Let $(X,\ms{O}_X)$ denote the inductive limit of the directed system of semialgebraic spaces $\{N_j \cap X | j\in J\}$. 
Then $(X,\ms{O}_X)$ has a natural structure of a locally semialgebraic space. 
With this locally semialgebraic space structure, $X$ is called a {\em locally semialgebraic subspace of $M$}.
\end{propdef}

The concept of a locally semialgebraic subspace provides a formal framework for precisely introducing the idea of a {\em semialgebraic subset} within a locally semialgebraic space.

\begin{propdef}\label{D:propdef2}
Let $(M,\ms{O}_M)$ be a locally semialgebraic space over $R$. 
A subset $X\subset M$ is called {\em semialgebraic} if $X$ is a locally semialgebraic subset of $M$ and the pair $(X,\ms{O}_X)$, where $\ms{O}_X$ is the restriction of $\ms{O}_M$ to $X$, is a semialgebraic space. 
If a locally semialgebraic subset $X\subseteq M$ is semialgebraic if and only if there exists a semialgebraic subset $Y\subseteq M$ such that 
$X\subset Y$.
\end{propdef}

\begin{example}
We proceed the notation of Example~\ref{E:U_r's}. 
In Figure~\ref{F:exceptional}, we depict a curve $Z \subset \R^2$, defined by the graph of the piecewise linear function $f:\R\to \R$, where 
\[
f(x) = \frac{2}{\pi} \sin^{-1} ( \sin (\pi x)),\quad x\in \R. 
\]

\begin{figure}[htp]
\begin{center}
\scalebox{0.7}{
\begin{tikzpicture}

\node at (.25,1) {1};
\node at (.25,-1) {-1};

\node at (1,-.25) {1};
\node at (2,-.25) {2};
\node at (3,-.25) {3};
\node at (-1,-.25) {-1};
\node at (-2,-.25) {-2};
\node at (-3,-.25) {-3};
\draw[ultra thick, color=blue]  (-5/2,-1) -- (-3/2,1) -- (-1/2,-1) -- (0,0) -- (1/2,1) -- (3/2,-1) -- (5/2,1)  ;
\draw[ultra thick, <-> ] (-6,0) -- (6,0) ;
\draw[ultra thick, <->] (0,-3) -- (0,3) ;

\draw[ultra thick, dashed, color=blue] (5/2,1) -- (7/2,-1) -- (4,0);
\draw[ultra thick, dashed, color=blue] (-5/2,-1) -- (-7/2,1) -- (-4,0);

\end{tikzpicture}    
}
\end{center}
\caption{A (connected) locally semialgebraic subset $Z\subseteq \R^2$.}
\label{F:exceptional}
\end{figure}
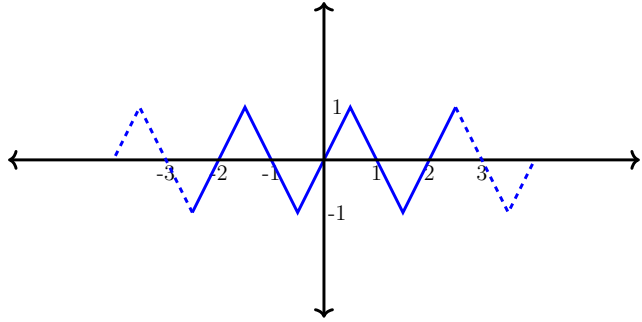
Clearly, for every $r\in \N$, the intersection $Z\cap U_r$ is a bounded semialgebraic subset of $U_r$. 
Indeed, $Z\cap U_r$ is a union of finitely many line segments. 
Hence, according to Remark~\ref{R:tocheckthat}, $Z$ is a locally semialgebraic subset of $(X,\ms{O}_X)$.
By a similar reasoning, we see that $Z$ is a locally semialgebraic subset of $(\R^2,\ms{O}_{\R^2})$ as well. 
However, we will show that $Z$ fails to qualify as a semialgebraic subset in either $(\R^2,\ms{O}_{\R^2})$ or $(X,\ms{O}_X)$. 

Assume towards a contradiction that $Z$ is a semialgebraic subset of $(\R^2,\ms{O}_{\R^2})$. 
Since intersection of two semialgebraic subsets is a semialgebraic subset, the intersection of $Z$ with the line $y=0$ in $\R^2$,
denoted by $Y$, must be a semialgebraic subset as well. 
But $Y$ consists of infinitely many isolated points in $\R^2$. 
Hence, $Y$ cannot be expressed by finitely many polynomial equalities or inequalities.
In particular, $Y$, hence $Z$, cannot a semialgebraic subset of  $(\R^2,\ms{O}_{\R^2})$.

The fact that $Z$ is not a semialgebraic subset of $(X,\ms{O}_X)$ follows from Definition~\ref{D:propdef2}.
Indeed, if $(Z,\ms{O}_X\vert_Z)$ were a semialgebraic space, then it would have a closed point. 
However, a basis for the generalized topology of $Z$ is given by $\{ Z\cap U_r:\ r\in \N \cup \{\infty\}\}$.
Since the complement of $Z\cap U_r$ ($r\in \N\cup \{\infty\}$) is either empty or unbounded, the generalized topology on $Z$ does not possess any closed points. 
This contradiction shows that $Z$ is not a semialgebraic subset of $(X,\ms{O}_X)$, also.
\end{example}

The preceding example generalizes when considering infinite locally semialgebraic spaces. In any such space $M$, there exists a locally semialgebraic subset $X \subsetneq M$ that does not qualify as semialgebraic. A proof of this fact mirrors the rationale employed in our previous example, but its details are omitted to avoid redundancy.

\begin{notation}
Let $X$ be a locally semialgebraic space. 
Recall that $\FT(X)$ denotes the family of all open locally semialgebraic subsets of $X$.
Recall also that by $\mathfrak{T}(X)$ we denoted the family of all locally semialgebraic subsets of $X$. 
Likewise, we will denote by $\mathfrak{S}(X)$ the family of all semialgebraic subsets of $X$. 
The family of all closed semialgebraic subsets of $X$ will be denoted by $\bar{\mathfrak{S}}(X)$.
Finally, let us remind that $\FS(X)$ stands for the family of all open semialgebraic subsets of $X$. 
\end{notation}

A (restricted) topological space $X$ is said to be {\em separated} if the image of the diagonal mapping 
\begin{align*}
\varDelta : X \longrightarrow X \times X
\end{align*} 
is a closed subset. 
In algebraic geometry, separated schemes hold significant prominence. 
Notably, every affine scheme is separated (see \cite[Proposition 4.1]{Hartshorne}). 
Furthermore, a widely accepted definition of an algebraic variety is that it is an integral separated scheme of finite type over an algebraically closed field
(see \cite[Example 4.0.2 and Proposition 4.10]{Hartshorne}).  
Serre proved in~\cite{GAGA} that a scheme defined over $\C$ is separated if and only if the associated analytic space is Hausdorff. 
These facts have direct consequences on a substantial family of semialgebraic spaces. 
Let $X$ be an algebraic variety defined over $R$. 
Let $M$ denote the set of $R$-rational points of $X$. 
Since all varieties are defined, at least locally, as the vanishing set of polynomials, $M$ has the structure of a semialgebraic space over $R$. 
Since $X$ is separated, $M$ is a separated topological space as well. 
(Further discussion on such semialgebraic spaces can be found in~\cite{DelfsKnebusch1981II}.)
For example, if $X$ denotes the complex affine line as an algebraic variety, then its set of $\R$-rational points is given by $M:=X(\R)= \R^2$. 
Then $M$ is separated with respect to its Zariski topology. 
But every Zariski closed subset of $M\times M$ is a closed subset in the restricted topology of $M\times M$. 
It follows that $M$ is separated with respect to its semialgebraic topology as well. 
For our first main theorem, we will utilize another separation property. 

\begin{definition}
A locally semialgebraic space $X$ is called {\em regular} if for every closed locally semialgebraic subset $A\subseteq X$ and every point $x\in X\setminus A$, there exist disjoint sets $U$ and $V$ from $\FT(X)$ such that $x\in U$, $A\subset V$.
\end{definition}

An important theorem of Robson~\cite{Robson} states that, a semialgebraic space $X$ is regular if and only if $X$ is affine.
In particular, for every $n\in \Z_+$, the affine space $(R^n,\ms{O}_{R^n})$ is a regular semialgebraic space.
However, a locally semialgebraic space need not be regular even if its underlying set is $R^n$ for some $n\in \Z_+$. 
\begin{example}
Let $X$ denote the locally semialgebraic space that we constructed in Example~\ref{E:U_r's}.
We will maintain the notation from that example. 
We claim that $X$ is not a regular locally semialgebraic space. 
Let $A$ denote the closed locally semialgebraic subset $X\setminus U_1$. 
Let $x$ denote the origin of $\R^2$. 
Then $x\in X\setminus A$. 
Given the nested nature of the $U_r$ sets, the set of all open subsets of $X$ is given by $\FT(X) = \{ U_r :\ r\in \N \} \cup \{ \R^2\}$. 
Let $r,s\in \Z_+$. 
Then we have $U_r \cap U_s = U_{\min \{r,s\}} \neq \emptyset$.
In other words, every two nonempty, open, locally semialgebraic subsets of $X$ intersect nontrivially.  
This means that we cannot separate $x$ from $A$ by using the nonempty open subsets of $X$.
Hence, $X$ is not regular.
\end{example}

Another topological notion that is naturally a part of semialgebraic geometry is the following. 

\begin{definition}
A locally semialgebraic space $X$ is called {\em paracompact} if it possesses an admissible covering $\{ X_i | i\in I\} \in \text{Cov}_X$ such that 
$\{ X_i | i\in I\}$ is a locally finite family, and every $X_i$ ($i\in I$) in this family is semialgebraic.
(Recall that we introduced the notion of a locally finite family in Definition~\ref{D:locallyfinite}.)
\end{definition}

Note that every semialgebraic space is paracompact. 

\medskip

Let $X$ be a locally semialgebraic space. 
A {\em family of supports on $X$} is a family $\Phi$ of closed subsets of $X$ such that 
\begin{enumerate}
\item[a)] every closed subset of a member of $\Phi$ is also an element of $\Phi$; 
\item[b)] $\Phi$ is closed under finite unions. 
\end{enumerate}
If, in addition, the following conditions are satisfied, then we call $\Phi$ a {\em paracompactifying family of supports}:
\begin{enumerate}
\item[c)] Every element of $\Phi$ is contained in some element of $\Phi \cap \mathfrak{T}(X)$. 
\item[d)] Every element of $\Phi \cap \mathfrak{T}(X)$ is regular and paracompact.
\item[e)] Every element of $\Phi$ has a neighborhood in $X$ which belongs to $\Phi$. 
\end{enumerate}

\begin{example}
Let $X$ be an affine semialgebraic space. 
Then the family $\Phi$ of all closed subsets of $X$ is a paracompactifying family of supports. 
\end{example}

\begin{remark}
The definition of a paracompactifying family of supports that we introduced above is the adaptation to the semialgebraic geometry of a similar notion from classical algebraic topology (see~\cite[Section 3.2]{Godement}). 
We will use the classical version in Section~\ref{S:Affinenash} when we focus on the sheaves on affine Nash manifolds. 
\end{remark}

Let $X$ be a locally semialgebraic space. 
Let $\mathscr{F}$ be an {\em abelian sheaf} on $X$, that is, a sheaf on $X$ with values in the category of abelian groups.
Let $s\in \Gamma(X,\mathscr{F})$. 
For every $x\in X$, the section $s$ defines an element in the stalk $\mathscr{F}_x$. 
The {\em support of $s$}, denoted $\text{supp }s$, is the set $\{ x\in X | s_x \neq 0\}$. 
Note that $\text{supp }s$ is a closed subset of $X$. 
For any family $\Phi$ of supports on $X$, we set 
\begin{align*}
\Gamma_{\Phi}(X,\mathscr{F}) := \{ s\in \Gamma (X, \mathscr{F}) | \ \text{supp }s\in \Phi \}.
\end{align*}
Then it is easy to check that the functor $\mathscr{F}\mapsto \Gamma_{\Phi}(X,\mathscr{F})$ is left exact. 
Hence, by using its right derived functors, 
$R^q \Gamma_\Phi (X, -)$, we obtain the {\em sheaf cohomology groups with supports in $\Phi$}.
In more mathematical notation, we will write 
\begin{align*}
H_\Phi^q ( X,\mathscr{F}) := R^q  \Gamma_\Phi (X,\mathscr{F})
\end{align*}
for the $q$-th cohomology group with supports in $\Phi$ of the sheaf $\mathscr{F}$. 
A sheaf $\mathscr{F}$ on $X$ is called $\Phi$-acyclic if $H_\Phi^q (X,\mathscr{F}) = 0$ for every $q\in \Z_+$.

Let $A$ be a closed subset of $X$. 
In our next definition, we will consider the restriction of a sheaf $\mathscr{F}$ onto $A$.
For convenience, let us recall its definition.
Then, the section ring $\Gamma(A,\mathscr{F})$ is defined as the inductive limit, 
\begin{align*}
\Gamma(A,\mathscr{F}):= \varinjlim_{\scriptsize{{A\subset U \in \ms{F}}}}\ms{F}_M(U).
\end{align*}

\begin{definition}
Let $\Phi$ be a paracompactifying family of supports on a locally semialgebraic space $X$. 
\begin{enumerate}
\item A sheaf $\mathscr{F}$ on $X$ is said to be {\em $\Phi$-soft} if the restriction map 
\begin{align*}
\Gamma(X,\mathscr{F})\longrightarrow \Gamma(A,\mathscr{F})
\end{align*} 
is surjective for every $A\in\Phi$. If $\Phi$ consists of all closed subsets of $X$, then $\mathscr{F}$ is called {\em soft}.
\item 
A sheaf $\mathscr{F}$ on $X$ is called {\em sa-flabby} if the restriction map 
\begin{align*}
\Gamma(X,\mathscr{F})\longrightarrow \Gamma(U,\mathscr{F})
\end{align*} 
is surjective for every $U\in\FT(X)$.
\end{enumerate}
\end{definition}

We now state several useful propositions from Delfs work.

\begin{proposition}\label{P:softisacyclic}
Let $X$ be a locally semialgebraic space. 
Let $\Phi$ be a paracompactifying family of supports on $X$. 
If $\mathscr{F}$ is $\Phi$-soft, then $\mathscr{F}$ is $\Phi$-acyclic.
\end{proposition}
\begin{proof}
This is a special case of the proof of~\cite[Proposition 4.17]{Delfs}.
\end{proof}

\begin{proposition}\label{P:softmoduleissoft}
Let $X$ be a locally semialgebraic space. 
Let $\Phi$ be a paracompactifying family of supports on $X$. 
Let $A$ be a principal ideal domain. 
Let $\mathscr{G}$ be a sheaf on $X$. 
If $\mathscr{F}$ is a flat and $\Phi$-soft sheaf of $A$-modules on $X$, then 
the sheaf $\mathscr{F}\otimes \mathscr{G}$ is $\Phi$-soft as well.
\end{proposition}
\begin{proof}
This is a restatement of~\cite[Proposition 4.22]{Delfs}.
\end{proof}

A locally semialgebraic space $X$ is said to be {\em Lindel\"of} if it possesses a countable admissible covering $\{ X_i | i\in I\} \in \text{Cov}_X(X)$.  
Clearly, every semialgebraic space is Lindel\"of.

\begin{proposition}\label{T:saflabbyisacyclic}
Let $X$ be a  Lindel\"of locally semialgebraic space. 
If $\mathscr{F}$ is an sa-flabby sheaf on $X$, then $\mathscr{F}$ is $\Phi$-acyclic for every locally semialgebraic family $\Phi$ of supports on $X$. 
\end{proposition}
\begin{proof}
The proof follows from~\cite[Theorem 4. 10]{Delfs}.
\end{proof}

We continue with the review of two more useful results about locally semialgebraic spaces. 
The first one is a partition of unity statement.
\begin{theorem}(\cite[Theorem 4.12]{DelfsKnebusch})
Let $X$ be a regular paracompact locally semialgebraic space. 
Let $\{ U_\lambda \to X | \lambda \in I \}$ be a locally finite covering of $X$ by open semialgebraic subsets. 
Then there exists a family $\{ \varphi_\lambda | \lambda \in I\}$ of locally semialgebraic functions 
$\varphi_\lambda : M\to [0,1]$ with $\text{supp}(\varphi_\lambda)\subset U_\lambda$ for every $\lambda \in I$ and 
\begin{align*}
\sum_{\lambda \in I} \varphi_\lambda (x)= 1
\end{align*}
for every $x\in X$.
\end{theorem}

The second result that we want to quote is a weak version of a locally semialgebraic analog of the well-known {\em Tietze's extension theorem} for continuous function.
Its proof heavily depends on the previous partition of unity theorem. 

\begin{theorem}(\cite[Theorem 4.13]{DelfsKnebusch})\label{T:Tietzeforlocallysa}
Let $X$ be a regular paracompact locally semialgebraic space. 
Let $f: A\to R$ be a locally semialgebraic function on a closed locally semialgebraic subset $A\subset M$.
Then there exists a locally semialgebraic function $g: M\to R$ such that $g\vert_A = f$. 
\end{theorem}

\begin{remark}\label{R:Tietzeforsa}
The special case of Theorem~\ref{T:Tietzeforlocallysa} for affine semialgebraic spaces is due to Efroymson~\cite{Efroymson} and Pecker~\cite{Pecker}. 
\end{remark}

We now proceed to introduce the group objects of the semialgebraic category. 
Let $G$ be a group. 
If $G$ has a locally semialgebraic space structure such that the group operations 
\begin{align*}
m: G\times G &\longrightarrow G \\
(x,y) &\longmapsto xy
\end{align*}
and 
\begin{align*}
i: G &\longrightarrow G \\
x&\longmapsto x^{-1}
\end{align*}
are morphisms of locally semialgebraic spaces, then we call $G$ a {\em locally semialgebraic group}.
A {\em semialgebraic group} is a locally semialgebraic group whose underlying locally semialgebraic space is a semialgebraic space. 
An {\em affine semialgebraic group} is a semialgebraic group whose underlying semialgebraic space is isomorphic to an affine semialgebraic space.
\begin{example}
Every real algebraic group is an affine semialgebraic group.
\end{example}
 
We know from the work~\cite{Pillay1988} of Pillay that every semialgebraic group has the unique structure of a Lie group but not conversely. 
Indeed, already the unit circle $\mathbb{S}^1$ can be endowed with infinitely many different semialgebraic group stuctures some of which are not even embeddable into any euclidean space. 
The first example of such a structure was found by Mazur. 
For a detailed discussion of the example of Mazur, see~\cite[Chapter IV]{Shiota}.

\section{Ringed $D$-Spaces and Vector Bundles}\label{S:Ringed}

Throughout this section, we fix the following data: 
\begin{longtable}{l c l}
$k$ &:& a field of characteristic zero, \\
$D$ &:&  a finite dimensional unital $k$-algebra,\\
$\varepsilon_0 : D \to k$ &:& a $k$-algebra homomorphism,\\
$J$ &:& the kernel of $\vep_0$ (which will be assumed to be nilpotent),\\ 
$X$ &:& a generalized topological space,\\
 $(X,\mathscr{C})$ &:& a locally ringed space over $k$.
\end{longtable}
As usual, if there is no danger for confusion, we will write $X$ instead of $(X,\mathscr{C})$. 
We are now ready to introduce {\em ringed $D$-spaces}.

\begin{definition}
A sheaf of $k$-algebras $\ms{A}$ on $X$ is called a {\em ringed $D$-space structure on $X$} if there exist a map of sheaves of $k$-algebras 
$\varepsilon : \ms{A} \rightarrow \mathscr{C}$ and a covering $\{ U_i  | i\in I \}$ of $X$ by open subsets, where, for each $i\in I$, there is a sheaf isomorphism 
\begin{align*}
\tau_i : \ms{A} \vert_{U_i} \longrightarrow (\mathscr{C}\otimes D) \vert_{U_i}
\end{align*}
such that $$(id \otimes \varepsilon_0) \circ \tau_i = \varepsilon \vert_{U_i}.$$
\end{definition}
\medskip

We will show that the notion of a ringed $D$-space is closely related to a notion of a vector bundle over $X$ with values in $D$. 
We proceed with the assumption that $X$ is connected in an appropriate sense in its Grothendieck topology. 
This assumption is not essential but simplifies our discussion. 
Also, temporarily, we consider only the vector space structure on $D$.

\begin{definition}\label{D:vbofclassC}
A {\em vector bundle over $X$ of class $\mathscr{C}$} is a 
sheaf $\ms{E}$ of $\ms{C}$-modules on $X$ for which there is an open covering $\{X_\alpha \to X | \alpha \in I\}$ satisfying the following conditions:
\begin{enumerate}
\item For every $\alpha \in I$, there is an isomorphism of sheaves of $\mathscr{C}$-modules, 
\begin{align}\label{A:vb1}
\tau_\alpha : \mathscr{E} \vert_{X_\alpha} \longrightarrow \mathscr{C} \vert_{X_\alpha} \otimes D.
\end{align}
\item For $\{ \alpha, \beta\}\subset I$, the map of sheaves of $\mathscr{C}$-modules, 
\begin{align}\label{A:vb2}
g_{\alpha \beta}:= \tau_\alpha \circ \tau_\beta^{-1} : \ms{C} \vert_{X_\alpha \cap X_\beta} \otimes D\longrightarrow \ms{C} \vert_{X_\alpha \cap X_\beta} \otimes D,
\end{align}
called the {\em transition function of $\ms{E}$ on $X_\alpha \cap X_\beta$}, is an isomorphism of $\ms{C}$-modules.
\end{enumerate}
In this case, the general linear group $\text{GL}(D)$ is said to be the {\em structure group of $\ms{E}$}.
\end{definition}

We follow the notation of Definition~\ref{D:vbofclassC}.
For every $x\in X_{\alpha} \cap X_{\beta}$, the transition function $g_{\alpha,\beta}$ in (\ref{A:vb2}) gives a $\ms{C}_x$-linear automorphism 
$(g_{\alpha,\beta})_x$ of $\ms{C}_x \otimes D$. 
Since $\ms{C}_x$ is a local $k$-algebra, and since $D$ is a $k$-vector space, such an automorphism is given by an element of the general linear group $\text{GL}(D)$. 
In fact, the data of transition functions of $\ms{E}$ is equivalent to a data of maps, 
\begin{align*}
\{ X_{\alpha}\cap X_{\beta}\to \text{GL}(D)\ |\ \{\alpha,\beta \} \subset I\},
\end{align*}
satisfying the {\em \v{C}ech cocycle conditions}, $(g_{\alpha,\alpha})_x = \text{id}_D$ and $(g_{\beta,\gamma})_x (g_{\alpha,\beta})_x
= (g_{\alpha,\gamma})_x$ for every $\{\alpha,\beta,\gamma \} \subset I$.

\begin{example}
In this example, we will find out that the ordinary topological vector bundles are among the examples of our Definition~\ref{D:vbofclassC}.
Let $X$ be a topological space. 
Let $k:=\R$.
Let $\ms{C}^0$ denote the sheaf of continuous $k$-valued functions on $X$. 
Let $\ms{E}$ be a locally free sheaf of $\ms{C}^0$-modules on $X$ such that there is an open cover $\{X_\alpha  | \alpha \in I\}$ of $X$ 
satisfying the properties in (\ref{A:vb1}) and (\ref{A:vb2}).
After fixing a basis for $D$, we may assume that $D$ is isomorphic to $\R^r$ for some $r\in \Z_+$. 
Then we have the isomorphism,
\begin{align*}
\ms{E}  \vert_{X_\alpha} \cong 
\underbrace{\ms{C}^0 \vert_{X_\alpha} \oplus \cdots \oplus \ms{C}^0 \vert_{X_\alpha}}_\text{$r$-copies} \qquad (\alpha \in I).
\end{align*}
Thus, locally, the transition functions are given by the elements of $\text{GL}_r(\R)$.
They satisfy the \v{C}ech cocycle conditions. 
Therefore, we recover the well-known definition of a topological vector bundle on $X$. 
\end{example}

The structure group of a vector bundle over $X$ of class $\ms{C}$ is a group object in the category that $X$ belongs to. 
This is forced upon us by the transition functions of the vector bundle.
In this sense, the structure group of our previous example, that is $\text{GL}_r(\R)$, is viewed as a topological group. 
Of course, $\text{GL}_r(\R)$ is a real algebraic group. 
Hence, it is also a locally semialgebraic group. 
In this regard, for us, a locally semialgebraic vector bundle on a locally semialgebraic space $X$ is a vector bundle $\ms{A}$ over $X$ of type $\ms{O}_X$ whose structure group is $\text{GL}_r(\R)$.
\medskip

We proceed with the assumption that $D$ is a $k$-algebra.
Let $\text{Aut}(D)$ denote the group of all $k$-algebra automorphisms of $D$. 
Then $\text{Aut}(D)$ is an algebraic subgroup of $\text{GL}(D)$.
Let $\ms{E}$ be a vector bundle over $X$ of class $\mathscr{C}$. 
In the extreme case that (the images in the stalks of) every transition function of $E$ as in (\ref{A:vb2}) is contained in $\text{Aut}(D)$, 
the sheaf $\ms{E}$ can be viewed as a sheaf of $\ms{C}$-algebras so that the map (\ref{A:vb1}) is an isomorphism of sheaves of $\ms{C}$-algebras. 
Let $\vep_0 : D\to k$ be a $k$-algebra homomorphism.
The automorphisms of $D$ that preserve $\vep_0$ form a subgroup of $\text{Aut}(D)$, denoted by $\text{Aut}(D,\vep_0)$.
Let us assume that $\ms{E}$ is a vector bundle over $X$ of class $\ms{C}$ such that all transition functions of $\ms{E}$ preserve
the map $\text{id}\otimes \vep_0$. 
\begin{definition}
If $\ms{E}$ is a vector bundle over $X$ of class $\ms{C}$ such that all transition functions of $\ms{E}$ preserve the map $\text{id}\otimes \vep_0$,
then we say that {\em $\ms{E}$ has the structure group $\text{Aut}(D,\vep_0)$}.
\end{definition}
Let $\ms{E}$ be a vector bundle over $X$ of class $\ms{C}$ with structure group $\text{Aut}(D,\vep_0)$. 
Then we have a morphism $\vep : \ms{E}\to \ms{C}$ such that locally on $U_\alpha$ it is given by 
\begin{align*}
\vep = (\text{id}\otimes \vep_0) \circ \tau_\alpha.  
\end{align*}
Let us put these definitions together in the form of a lemma.

\begin{lemma}\label{L:vbisaDspace}
Let $\ms{E}$ be a vector bundle over $X$ of class $\ms{C}$.
If all transition functions of $\ms{E}$ are contained in $\text{Aut}(D,\vep_0)$, 
then $\ms{E}$ has the structure of a ringed $D$-space on $X$.

Conversely, a ringed $D$-space $\ms{E}$ on $X$ is a vector bundle of class $\ms{C}$ if $\ms{E}$ is locally free as a sheaf of $\ms{C}$-modules. 
\end{lemma}

Let $\ms{A}$ be a sheaf of $k$-algebras on $X$. 
Let $\{U_\alpha | \alpha \in I\}$ be an open cover of $X$. 
Then the assignment $U_\alpha \mapsto Z(\ms{A}(U_\alpha))$, where $Z(\cdot)$ stands for the center, gives us a presheaf of $k$-algebras on $X$.
We denote the sheafification of this presheaf by $Z(\ms{A})$. 
We note in passing that, in general, the existence of the sheafification of a presheaf with values in the category of $k$-algebras is a nontrivial fact;
see~\cite[\href{https://stacks.math.columbia.edu/tag/00YR}{Tag 00YR}]{stacks-project}.
In our next proposition, a criterion for checking when a ringed $D$-space is locally free is presented.

\begin{proposition}\label{P:BRP1}
Let $\vep : \ms{A}\to \ms{C}$ be a ringed $D$-space structure on $X$. 
If $\varphi : \ms{C} \to Z(\ms{A})$ is a map of sheaves of $k$-algebras such that $\vep \circ \varphi = \text{id}_{\ms{C}}$, 
then $\ms{A}$ is a locally free sheaf of $\ms{C}$-modules. 
In this case, the structure group of $\ms{A}$ is contained in $\text{Aut}(D,\vep_0)$.
\end{proposition}

The proof of~\cite[Proposition 1]{BlattnerRawnsley} relies on the following useful lemma.

\begin{lemma}(\cite[Lemma 1]{BlattnerRawnsley})\label{L:Lemma1}
Let $\theta : D\to D$ be a unipotent endomorphism of $D$. 
Then $\log \theta := \sum_{p=1}^\infty (-1)^{p-1} (\theta - \text{id})^p / p$ is a nilpotent derivation of $D$.
If $\delta$ is a nilpotent derivation of $D$, then $\exp \delta := \text{id} + \sum_{p=1}^\infty \delta^p /p!$ is a unipotent automorphism of $D$.
The operations $\log$ and $\exp$ are inverses of each other. 
Furthermore, $\log$ is a bijection between the set of unipotent automorphisms of $D$ and the set of nilpotent derivations of $D$. 
\end{lemma}

\begin{proof}[Proof of Proposition~\ref{P:BRP1}]
This proof is identical to the proof of~\cite[Proposition 1]{BlattnerRawnsley}.
We record here since it is instructive.

Let us define an endomorphism $\theta_\alpha: \ms{C} \vert_{U_\alpha} \otimes D \to \ms{C} \vert_{U_\alpha} \otimes D$ by setting 
locally $\theta_\alpha (c\otimes d):=\tau_\alpha (\varphi(c))1\otimes d$, where $\tau_\alpha$ is as in (\ref{A:vb1}).
Then we have 
\begin{align*}
\theta_\alpha (c\otimes d) -  c\otimes d = (\tau_\alpha(\varphi(c))  - c\otimes 1) \cdot (1\otimes d).
\end{align*}
Notice that 
\begin{align*}
\text{id}\otimes \vep_0  (\tau_\alpha(\varphi(c)) - c\otimes 1) = (\vep \circ \varphi) (c)  - c = 0.
\end{align*}
In other words, $\tau_\alpha(\varphi(c)) - c\otimes 1$ is an element of $\ms{C}\vert_{U_\alpha}\otimes J$.
It follows that $\theta_\alpha - \text{id}$ maps $\ms{C}\vert_{U_\alpha}\otimes J^p$ into $\ms{C}\vert_{U_\alpha}\otimes J^{p+1}$.
Hence, we have $(\theta_\alpha -\text{id})^{n+1} = 0 $ for some $n\in \Z_+$. 
This means that $\theta_\alpha$ is a unipotent endomorphism of sheaves of algebras.
Then by Lemma~\ref{L:Lemma1}, $\theta_\alpha$ is of the form $\exp \delta_\alpha$ for some 
nilpotent endomorphism $\delta_\alpha$. Hence, $\theta_\alpha$ is a sheaf automorphism. 
Let $\sigma_\alpha : \ms{A}\vert_{U_\alpha} \to \mc{C}\vert_{U_\alpha} \otimes D$ denote the map defined by 
\begin{align*}
\sigma_\alpha := \theta^{-1}_\alpha \circ \tau_\alpha.
\end{align*}
It is easy to see that $\sigma_\alpha$ is an isomorphism of sheaves of $\ms{C}$-modules. 
Thus, $\ms{A}$ is locally free with transition functions $\sigma_\alpha \circ \sigma_\beta^{-1}$.
These transition functions have values in $\text{Aut}(D)$.
The following calculation shows that $\theta_\alpha^{-1}$, hence, $\theta_\alpha$ preserves $\text{id}\otimes \vep_0$:
\begin{align*}
\text{id}\otimes \vep_0 (\theta_\alpha (c\otimes d)) &= \text{id}\otimes \vep_0 (\tau_\alpha (\varphi(c)) 1\otimes d)\\
&= (\vep \circ \varphi) (c) \vep_0(d) \\
&=c \vep_0(d) \\
&= \text{id}\otimes \vep_0 (c\otimes d).
\end{align*}
It follows that $\sigma_\alpha \circ \sigma_\beta^{-1}$ has in fact values in $\text{Aut}(D,\vep_0)$. 
This finishes the proof of the proposition.
\end{proof}

We see from Proposition~\ref{P:BRP1} that, to show that a ringed $D$-space $\ms{A}$ is a vector bundle, 
it suffices to show that the map $\vep : Z(\ms{A}) \to \ms{C}$ is split as a surjective map of sheaves of $\ms{C}$-algebras. 
Indeed, this is the main strategy of the article of Blattner and Rawnsley,~\cite{BlattnerRawnsley} for proving the Batchelor's theorem in a general setup. 
We will adapt their strategy into our context. 
\medskip

Recall our notation that $J$ is the kernel of a $k$-algebra homomorphism, $\vep_0: D\to k$. 
We fix some additional notation. 
First, we introduce the notation corresponding to the ideal $J$ in the center of $D$, that is, 
\begin{align*}
J_0 := J \cap Z(D).
\end{align*}
If $\vep : \ms{A}\to \ms{C}$ is a ringed $D$-space structure on $X$, then we denote by $\ms{I}$ the ideal sheaf 
\begin{align*}
\ms{I} := \ker \vep.
\end{align*}
Then we have an ideal sheaf defined by 
\begin{align*}
\ms{I}_0 := \ms{I} \cap Z(\ms{A}).
\end{align*}
Finally, for $p\in \Z_+$, we denote by $\text{Der}(\ms{C}, \ms{I}_0^p / \ms{I}_0^{p+1}))$ the sheafification of the abelian presheaf of ``derivations of $\ms{C}$ into $\ms{I}_0^p / \ms{I}_0^{p+1}$.'' Here, by a {\em derivation of $\ms{C}$ into $\ms{I}_0^p / \ms{I}_0^{p+1}$}, we mean a map 
$\ms{C}\to \ms{I}_0^p / \ms{I}_0^{p+1}$
of sheaves of $\ms{C}$-algebras that has, locally, the well-known Leibniz property. 
Expressed using this notation, the principal theorem from \cite{BlattnerRawnsley} takes the form of the following assertion for the ringed $D$-spaces associated with a generalized topological space.

\begin{theorem}\label{T:maintool}
We maintain our previous data and notation. 
In addition, we assume that $J$ is a nilpotent ideal. 
If $\ms{A}$ is a ringed $D$-space structure on $X$ such that for every $p\in \Z_+$ the cohomology group  
\begin{align*}
H^1 (X, \text{Der}(\ms{C}, \ms{I}_0^p / \ms{I}_0^{p+1}))
\end{align*} 
is trivial, then there is a surjective map $\varphi : \ms{C} \to Z(\ms{A})$ of sheaves of algebras on $X$ such that $\vep \circ \varphi = \text{id}_{\ms{C}}$. 
In this case, as a sheaf of $\ms{C}$-modules, $\ms{A}$ is a vector bundle over $X$ of type $\ms{C}$ with structure group in $\text{Aut}(D,\vep_0)$. 
\end{theorem}

The proof of Theorem~\ref{T:maintool} develops in the same way as the proof of the corresponding theorem of~\cite{BlattnerRawnsley} does. 
The only difference is that we consider generalized topologies instead of topological spaces. 
To avoid an excessive amount of repetition, for the details of the proof, we refer the reader to the article~\cite{BlattnerRawnsley}.

\section{Proof of Theorem~\ref{T1:intro}}\label{S:Locally}

In this section we will provide a proof of our first result, Theorem~\ref{T1:intro}.

\medskip

Racall that a $\Z_2$-graded commutative $R$-algebra is a $\Z_2$-graded $R$-algebra $A$ such that 
for every $a\in A_m$ and $b\in A_n$, where $\{m,n\}\subseteq \{\overline{0},\overline{1}\}$, the following anti-commutation relation holds:
\begin{align}\label{A:anticommutation}
ab = (-1)^{m+n} ba. 
\end{align}
For our purposes, the most important example is the full exterior algebra $\LL R^s$. 
Let us denote it by $D$. 
Then we have 
\begin{align*}
D:= R^s\  \oplus \ \LL^1 R^s \ \oplus \ \LL^2 R^s \ \oplus \cdots \oplus \ \LL^s R^s. 
\end{align*}
The $\Z_2$-grading on $D$ is given by the decomposition $D = D_{\overline{0}} \oplus D_{\overline{1}}$,
where 
\begin{align*}
D_{\overline{0}}:= (\LL  R^s)_{\overline{0}} = \bigoplus_{i \ : \ \text{even} } \LL R^i\qquad \text{and}  \qquad
D_{\overline{1}}:= (\LL  R^s)_{\overline{1}} = \bigoplus_{i \ : \ \text{odd} } \LL R^i. 
\end{align*}
It is easy to check that, under this $\Z_2$-grading, the anti-commutation rule (\ref{A:anticommutation}) holds for $D$. 
In other words, $D$ is a $\Z_2$-graded commutative $R$-algebra. 
At the same time, it is also easy to check from (\ref{A:anticommutation}) that the center of $D$ is given by $D_{\overline{0}}$.

Let $e_1,\dots, e_s$ denote the standard basis vectors for $R^s$.
Then as an $R$-algebra, $D$ is generated by $1,e_1,\dots, e_s$. 
In this notation, there is a unique $R$-algebra homomorphism $\vep_0 : D\to R$ such that $\vep_0 ( e_i ) = 0$ and 
$\vep_0(1) = 1$. Let $J$ denote the kernel, $\ker \vep_0$.
Evidently, $J$ is a nilpotent ideal of $D$.

We proceed to show that the hypothesis of Theorem~\ref{T:maintool} is satisfied for these choices.

\begin{lemma}\label{L:H1vanishesforsa}
Let $\vep : \ms{A}\to \ms{O}_X$ be a ringed $D$-space structure on a regular, paracompact locally semialgebraic space $X$. 
Let $\ms{I}$ denote the ideal sheaf defined by 
\begin{align*}
\ms{I} := \ker \vep.
\end{align*}
Let $\ms{I}_0$ denote the corresponding ideal sheaf in $Z(\ms{A})$, 
\begin{align*}
\ms{I}_0 := \ms{I} \cap Z(\ms{A}).
\end{align*}
Then for every $p\in \Z_+$, the first cohomology of $\text{Der}(\ms{O}_X, \ms{I}_0^p / \ms{I}_0^{p+1})$, that is, 
\begin{align*}
H^1 (X, \text{Der}(\ms{O}_X, \ms{I}_0^p / \ms{I}_0^{p+1}))
\end{align*} 
is trivial. 
\end{lemma}

\begin{proof}
We already know from Theorem~\ref{T:Tietzeforlocallysa} that the Tietze's extension theorem holds for our locally semialgebraic spaces.
This theorem readily implies that the sheaf $\ms{O}_X$ is a soft sheaf. 
Note that, since $\ms{O}_X$ is a sheaf of $R$-algebras, it has a natural structure of a sheaf of $\Z$-modules, where 
we view the elements of $\Z$ as (locally) constant functions on $X$.  
Since the sheaf $\text{Der}(\ms{O}_X, \ms{I}_0^p / \ms{I}_0^{p+1})$ is a sheaf of $\ms{O}_X$-modules, 
we see by Proposition~\ref{P:softmoduleissoft} that $\text{Der}(\ms{O}_X, \ms{I}_0^p / \ms{I}_0^{p+1})$ is a soft sheaf on $X$. 
In particular, it is an acyclic sheaf. 
This finishes the proof of our assertion. 
\end{proof}

\medskip

We are now ready to finish the proof of our Theorem~\ref{T1:intro}. 
Let us recall its statement for convenience: 
\medskip

Let $X$ be a regular, paracompact locally semialgebraic space. 
If $(X,\ms{A})$ is a $\Z_2$-graded semialgebraic space of odd-dimension $s$ on $X$, 
then there exists a locally semialgebraic vector bundle $E$ on $X$ such that $\ms{A}(X)$ is isomorphic, 
as a $\Z_2$-graded commutative algebra, to $\Gamma(X,\LL E)$.

\begin{proof}[Proof of Theorem~\ref{T1:intro}]
Recall that $D$ is the $\Z_2$-graded commutative algebra $\LL R^s$.  
Let us denote by $\text{Aut}_0(D,\vep_0)$ the subgroup of $\text{Aut}(D,\vep_0)$ consisting of its $\Z_2$-grading preserving elements. 
Since $\text{Aut}_0(D,\vep_0)$ is an algebraic subgroup of $\text{Aut}(D,\vep_0)$, which is defined over $R$, 
it is also an affine semialgebraic subgroup. 
Now, following the general idea of the proof of the Batchelor's theorem in~\cite[Remark 2]{BlattnerRawnsley}, 
we use the cohomology vanishing result. 
To this end, we apply Lemma~\ref{L:H1vanishesforsa}.
Then Theorem~\ref{T:maintool} shows that our sheaf of $\Z_2$-graded commutative algebras $\ms{A}$ on $X$ can be viewed as the sheaf of 
sections of a vector bundle over $X$ of type $\ms{O}_X$ with structure group in $\text{Aut}_0(D,\vep_0)$.
To finish off the proof, we still need to show that the structure group of $\ms{A}$ can be reduced to the $\Z$-grading preserving subgroup of 
$\text{Aut}(D,\vep_0)$, which is $\text{GL}_s(R)$.
Of course, this is also an affine semialgebraic subgroup.
Also, the reasoning presented in the second half of~\cite[Remark 2]{BlattnerRawnsley} applies verbatim to the semialgebraic groups 
$\text{Aut}_0(D,\vep_0), \text{Aut}(D,\vep_0)$, and $\text{GL}_s(R)$.
Indeed, the quotients $\text{Aut}(D,\vep_0)/\text{GL}_s(R)$ and $\text{Aut}_0(D,\vep_0)/\text{GL}_s(R)$ are affine spaces.
Therefore, the former quotient can be be contracted by a semialgebraic homotopy to the latter quotient. 
Thus, we see that the structure group of $\ms{A}$ reduces to $\text{GL}_s(R)$ in the semialgebraic category as well. 
This finishes the proof of our first main theorem. 
\end{proof}

\section{Nash Manifolds and Their $\mathfrak{N}$-topologies}\label{S:Affinenash}

In this section, we prove our second main result, Theorem~\ref{T2:intro}. 
Towards this end, we will review, in more detail, the theory of affine Nash manifolds. 
Central to our development in this section is a class of functions on the affine Nash manifolds that behaves like the class of polynomial functions on an affine algebraic variety. 
We begin with a discussion of the relevant functions for semialgebraic spaces.
\medskip
 
Recall that previously we defined the morphisms of locally semialgebraic spaces as the morphisms of the underlying locally ringed spaces. 
Let us make this definition more explicit in the case of affine semialgebraic spaces. 

\begin{definition}
Let $n\in \Z_+$. 
Let $M$ be a semialgebraic subset of $\R^n$. 
If, in addition, $M$ is a smooth submanifold of $\R^n$, then $M$ is called a {\em Nash submanifold of $\R^n$}.
\end{definition}

Let $A$ (resp. $B$) be a semialgebraic subset in $\R^m$ (resp. in $\R^n$).
A map $f: A\to B$ is called a {\em semialgebraic map} if the graph of $f$ is a semialgebraic subset in $\R^m\times \R^n$.
Let $f: A\to B$ be a semialgebraic map, as defined above. Also, we assume that $A$ and $B$ are open sets. 
If, in addition, $f$ is a smooth map, then $f$ is called a {\em Nash map}. 
A {\em Nash diffeomorphism} is a bijective Nash map $f: A\to B$ whose inverse $f^{-1}:B\to A$ is a Nash map. 
Notice that these definitions regarding Nash maps did not use Nash submanifolds of $\R^n$. 
To clarify this picture, let us reintroduce the Nash submanifolds of $\R^n$ by using Nash maps. 

A semialgebraic subset $M\subset \R^n$ is said to be a Nash submanifold of $\R^n$ of dimension $d$ if, for every $x\in M$, there exists a 
Nash diffeomorphism $\varphi$ from an open semialgebraic neighborhood $\varOmega$ of the origin $0\in \R^n$ onto an open semialgebraic neighborhood $\varOmega'$ of $x$ in $\R^n$ such that $\varphi (0)=x$ and 
\[
\varphi (\R^d \times \{0\} \cap \varOmega) = M\cap \varOmega'.
\]
Clearly, every open semialgebraic subset $X$ of a Nash submanifold of $\R^n$ of dimension $d$ is a Nash submanifold of $\R^n$ of dimension $d$.

\begin{remark}\label{R:everyNashisclosed}
We would like to highlight a characteristic of affine Nash manifolds that resembles a property found in every affine algebraic variety defined over an algebraically closed field.
It is well-known that an affine algebraic variety $X\subset \C^n$ need not be closed in $\C^n$. 
For example, the set $X:=\C \setminus \{0\}$ is a Zariski open subset of $\C$. 
But $X$ is isomorphic to the ``hyperbola'' $V:=\{ (a,b)\in \C^2 :\ ab=1\}$ in $\C^2$ via the natural morphism 
\begin{align*}
X &\longrightarrow V,\\ 
x &\longmapsto (x,1/x).
\end{align*}
In other words, although $X$ is not a closed subset of $\C$, it is isomorphic to a Zariski closed subvariety in a higher dimensional affine space. 
Likewise, a Nash submanifold of $\R^n$ ($n\in \Z_+$) need not be closed in $\R^n$.
Nonetheless, every Nash submanifold of $\R^n$ is Nash diffeomorphic to a closed Nash submanifold of some $\R^m$ ($m\in \Z_+$).
A proof of this fact is given in~\cite[Corollary I.4.3]{Shiota}. 
We note also that, since every open semialgebraic subset $X$ of a Nash submanifold $M\subseteq \R^n$ has the structure of a Nash submanifold of $\R^n$, $X$ is isomorphic to a closed Nash submanifold of some $\R^m$ ($m\in \Z_+$).
\end{remark}
 
Let $M$ be a Nash submanifold of $\R^n$ for some $n$. 
A {\em Nash function} on an open semialgebraic subset $U\subseteq M$ is a smooth semialgebraic function. 
In fact, by a theorem of Malgrange~\cite{Malgrange}, every such map is a real analytic map. 
Thus, our definition of a Nash function here agrees with the one that is given in the introduction section. 
The Nash functions defined on an open semialgebraic subset $U\subseteq M$, denoted by $\ms{N}(U)$, form a subring of the ring of all functions defined on $U$. 
In fact, $\ms{N}(U)$ is a noetherian ring. 
This fact was originally proved independently by Efroymson in~\cite{Efroymson1974} and by Riseler in~\cite{Risler} for globally defined Nash functions on Nash submanifolds of $\R^n$. It follows from our Remark~\ref{R:everyNashisclosed} that it holds for every open semialgebraic subset of a Nash submanifold $M$ of $\R^n$. 
Also, it is not difficult to check that the assignment $U\mapsto \ms{N}(U)$, $U\in \FS(M)$ defines a sheaf on $M$, denoted by $\ms{N}_M$. 
We call it the {\em sheaf of Nash functions on $M$.}

\begin{definition}
Let $M$ be a Nash submanifold of $\R^n$.
A {\em Nash subset of $M$} is a semialgebraic subset of the form 
\begin{align*}
\mathcal{Z}_M(f_1,\dots, f_r) := \{ x\in M | f_1(x) = \cdots = f_r(x) = 0\},
\end{align*}
where $f_1,\dots, f_r$ are Nash functions from $M$ to $\R$. 
\end{definition}

We continue with the assumption that $M$ is a Nash submanifold of $\R^n$ for some $n$.
Since $\ms{N}(M)$ is a noetherian ring, every Nash subset of $M$ is defined by an ideal of $\ms{N}(M)$. 
Recall from the introduction section that the $\mathfrak{N}$-topology on $M$ is the restricted topology whose closed subsets are the Nash subsets of $M$. 
Hence, we know that every $\mathfrak{N}$-closed subset of $M$ is a zero-set of the form 
\begin{align*}
\mathcal{Z}_M(I) := \{ x\in M | f(x)= 0\text{ for every $f\in I$}\},
\end{align*}
where $I$ is an ideal of $\ms{N}(M)$. 

\medskip

Since the intersection of a family of closed sets in a restricted topology may not be a closed set in that restricted topology, the closure operation may not exist. 
Fortunately, in the restricted topology of semialgebraic subsets, every semialgebraic subset has a semialgebraic closure. 
A proof of this fact can be found in~\cite[Corollary 2.2.7]{AizenbudGourevitch2008}. 
Of course, in the case of the $\mathfrak{N}$-topology, we are not concerned about this remark since $\mathfrak{N}$-topology is indeed a topology.

\begin{notation}\label{N:restrictions}
The set of all open (resp. closed) subsets of the $\mathfrak{N}$-topology on $M$ will be denoted by $\FN(M)$ (resp. by $\bar{\mathfrak{N}}(M)$). 
By abuse of notation, we denote the sheaf obtained by restricting the sheaf of Nash functions $\ms{N}_M$ to the $\FN(M)$-topology by $\ms{N}_M$. 
Likewise, we denote the sheaf obtained by restricting the sheaf of continuous semialgebraic functions $\ms{O}_M$ to the $\FN(M)$-topology by $\ms{O}_M$ as well. 
\end{notation}

Since every Nash function is a continuous semialgebraic function, every $\mathfrak{N}$-closed subset of $M$ is an $\mathfrak{S}$-closed subset as well. 
However, as we demonstrate it in the following example, the containment $\bar{\mathfrak{N}}(M) \subset \bar{\mathfrak{S}}(M)$ can be a strict inclusion. 

\begin{example}
Let $M$ denote $\R^1$ with its unique affine Nash manifold structure. 
Let $U$ denote the interval $(0,\infty)$ in $\R^1$.
The semialgebraic closure of $U$ is given by the half-closed interval, 
\begin{align*}
[0,\infty) := \{ x \in \R^1 | \ x \geq 0 \}.
\end{align*}
We claim that $[0,\infty)$ cannot be an $\mathfrak{N}$-closed subset of $M$. 
Let us assume towards a contradiction that $[0,\infty) \in \bar{\mathfrak{N}}(M)$ holds. 
Then we find finitely many Nash functions $f_1,\dots, f_r$ on $M$ such that $[0,\infty) = \bigcap_{i=1}^r f^{-1}_i(0)$. 
In particular, since $[0,\infty)$ is an unbounded closed subset of $M$, there is a nonempty open subinterval $I \in \FN(M)$ such that 
$I \subset f_i^{-1}(0)$ for each $i\in \{1,\dots, r\}$.
But since $M$ is semialebraically connected, we know the following implication from~\cite[Proposition 8.1.9]{BochnakCosteRoy}: 
\begin{align*}
f_i \vert_{I} = 0\implies f_i = 0 \text{ on $M$ for every $i\in \{1,\dots, r\}$}.
\end{align*}
This means that $\bigcap_{i=1}^r f^{-1}_i(0) = M$, which is absurd. 
Therefore, the $\mathfrak{N}$-topology and the $\mathfrak{S}$-topology on $M$ are different topologies. 
\end{example}

\medskip

In~\cite[Definition 3.2.4]{AizenbudGourevitch2008}, an {\em $\R$-space} is defined as a ringed space $(M,\ms{O}_M)$, where $M$ is a restricted topological space, and $\ms{O}_M$ is a sheaf of $\R$-algebras over $M$ which is a subsheaf of the sheaf of all continuous real-valued functions on $M$. 
In~\cite[Example 3.1.1]{AizenbudGourevitch2008}, a ringed space of the form $(M,\ms{N}_M)$, where $M$ is a Nash submanifold of $\R^n$ ($n\in \Z_+$) and $\ms{N}_M$ is the sheaf of Nash functions defined on the restricted topology of $M$ is called 
the {\em $\R$-space} of the Nash submanifold $M$ of $\R^n$. 
Then an {\em affine Nash manifold} is defined as an $\R$-space which is isomorphic to the $\R$-space of a Nash submanifold of $\R^n$ for some $n\in \Z_+$. 

\begin{remark}
In the preliminaries section we introduced the affine semialgebraic groups as the group objects in the category of affine semialgebraic spaces. 
An {\em affine Nash group} can be defined similarly, as a group object in the category of affine Nash manifolds. 
\end{remark}

\begin{definition}
Let $M$ be a Nash submanifold of $\R^n$. 
We call $(M,\ms{N}_M)$ the {\em $\R$-space of $M$ in the $\mathfrak{N}$-topology}.
An {\em affine Nash manifold} is a ringed space that is isomorphic to the 
{\em $\R$-space in the $\mathfrak{N}$-topology} of a Nash submanifold of $\R^n$ for some $n\in \Z_+$. 
\end{definition}

\begin{example}
Every smooth real algebraic variety is an affine Nash manifold. 
The curve $C$ in $\R^2$ defined by the polynomial equation $y^2=x^3$ is an affine semialgebraic space but not an affine Nash manifold since it has a singularity at the origin. 
\end{example}

We proceed with a review of some useful facts about the semialgebraic coverings of $M$.

Let $F$ be a continuous semialgebraic function on $M$. 
Following~\cite[Notation A.3.1]{AizenbudGourevitch2008} we denote by $M_F$ the set $\{ x\in M:\ F(x)\neq 0\}$.
If the restriction of $F$ to $M_F$ is a positive Nash function, then $F$ is said to be a {\em basic semialgebraic function}.
Let $\{F_1,\dots, F_m\}$ be a set of basic semialgebraic functions on $M$ such that, for every $i\in \{1,\dots, m\}$,
$F_i$ is continuous (w.r.t. the $\mathfrak{S}(M)$-topology).
We call $\{F_1,\dots, F_m\}$ a {\em basic collection} if for every $x\in M$ at least one of the functions $F_1,\dots, F_m$ is greater than 1 at $x$. 
An open cover $M=\bigcup_{j=1}^m V_j$ is said to be a {\em refinement} of another open cover $M=\bigcup_{i=1}^n U_i$ 
if for any $j\in \{1,\dots, m\}$ there exists $i\in \{1,\dots,n\}$ such that $V_j \subset U_i$. 
We now state a result of Aizenbud and Gourevitch as a lemma so that we can refer to it easily. 
\begin{lemma}(\cite[Proposition A.3.4]{AizenbudGourevitch2008})\label{L:shrink}
If $M=\bigcup_{i=1}^n U_i$ is a finite open semialgebraic cover of a Nash manifold $M$, then there exists a basic collection $\{F_1,\dots, F_m\}$ such that 
$M=\bigcup_{i=j}^m M_{F_j}$ is a refinement of $M=\bigcup_{i=1}^n U_i$.
\end{lemma}
We will apply this lemma to the stalks of the sheaf $\ms{N}_M$ with respect to the $\mathfrak{N}$-topology. 
The ring of sections of $\ms{N}_M$ on a closed subset $V\subset M$ is given by the ring of global sections of the restricted sheaf,
$i^{-1} \ms{N}_M$, where $i: V\hookrightarrow M$ is the closed embedding of $V$ into $M$. 
Equivalently, it is defined by the direct limit, 
\begin{align}\label{A:directlimit}
\ms{N}_M(V) := \varinjlim_{\scriptsize{{V\subset U \in \FN(M)}}}\ms{N}_M(U).
\end{align}
This means that, for $s\in \ms{N}_M(V)$, there exists a pair $(U,\tilde{s})$, where $U$ is an $\mathfrak{N}$-open set such that $V\subset U$, and $\tilde{s}$ is an element of $\ms{N}_M(U)$ such that $s(x)=\tilde{s}(x)$ for every $x\in V$.
Two such pairs $(U,\tilde{s})$ and $(U',\tilde{s}')$ represent the same section $s\in \ms{N}_M(V)$ if and only if there exists a third $\mathfrak{N}$-open set $U''$ such that $V\subset U'' \subseteq U\cap U'$ and $\tilde{s}(x) = \tilde{s}'(x)$ for every $x\in U''$.
By specializing (\ref{A:directlimit}) to the closed sets of the form, $V=\{x\}$, where $x$ is a point in $M$, we recover the definition of the {\em stalk of $\ms{N}_M$ at $x$}, which we denote by $\ms{N}_{M,x}$.

\begin{proposition}\label{P:locallyringed}
Let $M$ be an affine Nash manifold. 
Then the pair $(M,\ms{N}_M)$ with respect to the $\mathfrak{N}$-topology of $M$ is a locally ringed space over $\R$.
\end{proposition} 

\begin{proof}
Since any open cover in the $\mathfrak{N}$-topology is an open cover in the $\mathfrak{S}$-topology on $M$, 
all sheaf axioms hold automatically for $\ms{N}_M$ with respect to the $\mathfrak{N}$-topology. 
It follows that $(M,\ms{N}_M)$ is indeed a ringed space over $\R$. 
We proceed to show that the stalk $\ms{N}_{M,x}$ is a local ring for every $x\in M$. 
Let $U$ be an open semialgebraic neighborhood of $x$ in $M$.  
We view $U$ as an open subset of a semialgebraic cover of $M$. 
Then Lemma~\ref{L:shrink} shows that $x$ is contained in an open subset $M_F\subseteq U$,
where $F$ is a basic semialgebraic function on $M$.  
This argument shows that every open neighborhood of $x$ can be shrinked to an element of $\FN(M)$. 
Hence, we have the following identification of direct limits:
\begin{align*}
\varinjlim_{\scriptsize{{V\subset U \in \FS(M)}}}\ms{N}_M(U) =\varinjlim_{\scriptsize{{V\subset U \in \FN(M)}}}\ms{N}_M(U).
\end{align*}
But it is well-known that~\cite[Corollary 8.1.6]{BochnakCosteRoy} the stalk 
$\varinjlim_{\scriptsize{{V\subset U \in \FS(M)}}}\ms{N}_M(U)$
is a local ring. This finishes the proof of our assertion.  
\end{proof}

\medskip

We will now review a significant theorem that is very useful for the cohomological analysis of Nash manifolds. 
This theorem is known as ``Efroymson's extension theorem.''
\begin{lemma}(\cite{Efroymson, Pecker})\label{L:EP}(Efroymson's extension theorem)
Let $M$ be an affine Nash manifold. 
Let $\varphi: M\to \R$ be a Nash function on $M$.
Let $g$ be a Nash function that is defined on a semialgebraic neighborhood $U$ of $\varphi^{-1}(0)$. 
Then there exists a Nash function $h : M\to \R$ such that $g-h \vert_U$ is given by the product of $\varphi \vert_U$ and a Nash function on $U$. 
In partigular, $g$ and $h$ coincide on the closed Nash set $\varphi^{-1}(0)$. 
\end{lemma}

We will use Efroymson's extension theorem in the proof of Theorem\ref{T:Nashsoft}.
We recall its statement for convenience.
\medskip

If $M$ is an affine Nash manifold, then the sheaf $\ms{N}_M$ is a soft sheaf with respect to the $\mathfrak{N}$-topology on $M$. 
\begin{proof}[Proof of Theorem~\ref{T:Nashsoft}]
Let $V$ be an $\mathfrak{N}$-closed subset of $M$. 
Then there are finitely many Nash functions $f_1,\dots, f_m$ on $M$ such that 
$V=\{x\in M |\ f_1(x) =\cdots = f_m(x) = 0\}$.
We define a Nash function $\varphi :M\to \R$ by the formula, $\varphi:=f_1^2+\cdots +f_m^2$. 
Clearly, for $x\in M$, we have 
\begin{align*}
\varphi(x) = 0 \iff f_i (x) =0 \text{ for all $i\in \{1,\dots,m\}$}.
\end{align*}
In other words, $V$ is the zero set of $\varphi$, that is, $V= \varphi^{-1}(0)$.

We are now ready to show that the following natural map
\begin{align*}
\text{res}_{V} : \Gamma( M, \ms{N}_M ) &\longrightarrow \Gamma( V, \ms{N}_M )\\
h &\longmapsto h\vert_V
\end{align*}
is surjective. 
To this end, let $s$ be a section of $\Gamma( V, \ms{N}_M )$. 
Then there exists an $\mathfrak{N}$-open neighborhood $V\subset U$ and a section $\tilde{s}\in \ms{N}_M(U)$ such that 
$\tilde{s}(x) = s(x)$ for every $x\in V$.  
By Lemma~\ref{L:EP}, we know that there exists a Nash function $h:M\to \R$ such that $h- \tilde{s}$ is a product of $\varphi\vert_U$ and a Nash function on $U$.
It follows that $h \vert_V = s$ on $V$. 
In other words, we have $\text{res}_V(h) = s$. 
Hence, we showed that the map $\text{res}_V$ is surjective. 
This finishes the proof of our theorem.
\end{proof}

Let $X$ be a topological space having a paracompactifying family of supports. 
As we mentioned it in the introductory section, every flasque (flabby) sheaf on $X$ is a soft sheaf.   
Unfortunately, our sheaf of Nash functions $\ms{N}_M$ with respect to the $\mathfrak{N}$-topology is not flabby unless $M$ is finite. 
\begin{example}
Let $M$ denote $\R^1$ with its standard Nash manifold structure. 
Let $U$ denote the $\mathfrak{N}$-open subset, $U:=M\setminus \{0\}$. 
It is easy to check that the rational map $f: U\to \R$, $x\mapsto 1/x$ gives a section of $\ms{N}_M$ on $U$. 
However, no globally defined Nash function $s: M\to \R$ has the property that $s|_U (x)= f(x)$ for every $x\in U$. 
In other words, the restriction map $\text{res}_{V} : \Gamma( M, \ms{N}_M ) \to \Gamma( U, \ms{N}_M )$ is not surjective. 
\end{example}
However, by utilizing standard techniques from algebraic geometry, we are able to establish the anticipated cohomological properties of our sheaves. To this end, we attach an affine scheme to every affine Nash manifold. 

\begin{definition}
Let $M$ be an affine Nash manifold. 
Let $A$ denote the ring of globally defined Nash functions on $M$, that is, $A:=\ms{N}_M(M)$. 
Then we define the {\em Nash scheme of $M$}, denoted $\text{Nsp}(M)$, as the affine scheme $\text{Spec } A$. 
\end{definition}

Since the $\R$-algebra of global Nash functions on an affine Nash manifold $M$ is a noetherian ring, the associated Nash scheme is a scheme of finite type over $\R$. In particular, the set of closed points is dense in $\text{Nsp}(M)$ (see~\cite[Ch II, Exercise 3.14]{Hartshorne}).
We will show that the closed points of $\text{Nsp}(M)$ are in bijection with the points of $M$.

\begin{lemma}\label{L:closedpoints}
Let $M$ be an affine Nash manifold. 
Let $A$ denote the ring $\ms{N}_M(M)$. 
If $x$ is a point of $M$ and $I$ is the ideal defined by $I:=\{f\in A:\ f(x) = 0\}$, then $I$ is a maximal ideal. 
Conversely, if $I$ is a maximal ideal of $A$,
and $V$ is the vanishing set defined by $V:=\{x\in M :\ f(x) = 0\ \text{for every $f\in I$}\}$, then $V=\{x\}$ for some $x\in M$. 
\end{lemma}
\begin{proof}
Let $B$ denote the stalk of $\ms{N}_M$ at $x$. 
Then we have a surjective ring homomorphism $\varphi: A\to B$ defined by sending a global Nash function $s$ to its germ at $x$. 
We know from Proposition~\ref{P:locallyringed} that $B$ is a local ring. 
It is easy to see that the unique maximal ideal $\tilde{I}$ of $B$ consists of pairs $(U,s)$, where $U$ is an open neighborhood (w.r.t. $\mathfrak{N}$-topology) of $x$ and $s\in \ms{N}_M(U)$ is such that $s(x) = 0$. 
Since $\varphi$ is surjective, $\varphi^{-1}(\tilde{I})$ is a maximal ideal of $A$. 
But $\varphi^{-1}(\tilde{I})$ consists of $f\in A$ such that $f(x) = 0$. 
In other words, we have $I = \varphi^{-1}(\tilde{I})$. 
Hence, the proof of our first assertion follows. 

Our second assertion is stated in~\cite[Corollary 8.6.3]{BochnakCosteRoy}.
Hence, the proof of our lemma is finished.
\end{proof}

\medskip

Lemma~\ref{L:closedpoints} shows that the Nash scheme of $M$ is packed with much information about the geometry of $M$.

\begin{lemma}\label{L:homeomorphism}
Let $M$ be an affine Nash manifold. Let $A:=\ms{N}_M(M)$. 
Let $\beta: M\to \text{Nsp}(M)$ be the map defined by $\beta(x) := I_x$,
where $I_x$ is the maximal ideal of $A$ corresponding to $x\in M$. 
Then $\beta$ is a homeomorphism onto the subspace of closed points of $\text{Nsp}(M)$.
\end{lemma}
\begin{proof}
Lemma~\ref{L:closedpoints} shows that $\beta$ is a bijection of $M$ onto the set of closed points of $\text{Nsp}(M)$. 
If $V$ is a closed subset of $\text{Nsp}(M)$, then there exists a prime ideal $P\subset A$ such that $V = \{Q:\ \text{$Q$ is a prime ideal of $A$ containing $P$}\}$.
Then we have 
\begin{align*}
\beta^{-1}(V) &= \{x\in M :\ P\subseteq I_x \}\\
&=\{x\in M :\ f(x) = 0\text{ for every $f\in P$}\}.
\end{align*}
Since $A$ is noetherian, $P$ is finitely generated. 
Hence, we can write $\beta^{-1}(V)$ as the vanishing set of a single Nash function by using the argument in the first paragraph of the proof of Theorem~\ref{T:Nashsoft}.
This means that $\beta^{-1}(V)$ is a closed subset of $M$. 
Hence, $\beta$ is a continuous map. 
Conversely, if $W$ is a Nash closed subset of $M$, then we can write it as the zero set of a single Nash function $f: M\to \R$. 
In particular, the image of $W$ under $\beta$ consists of all maximal ideals $Q$ of $A$ such that $f\in Q$. 
Since the set of all prime ideals of $A$ that contains $f$ is a closed subset, the intersection of this closed subset with $\beta(M)$ is precisely the subset 
$t(W)$. Hence, $t(W)$ is a closed subset of the subspace of closed points of $\text{Nsp}(M)$.
This finishes the proof of the fact that $\beta$ is a bijective map such that both $\beta$ and its inverse are continuous. 
Therefore, $\beta$ is a homeomorphism. 
\end{proof}

As a consequence of the proof of Lemma~\ref{L:homeomorphism}, we note that the structure sheaf of the Nash scheme of an affine Nash manifold $M$ restricts to give the structure sheaf of $M$. We will make this more precise in our next theorem.
First, we recall the definition of a quasi-coherent sheaf of $\ms{N}_M$-modules. 

Let $V$ be an $\mathfrak{N}$-closed subset of an affine Nash manifold $M$. 
Let $A$ denote the ring of Nash functions, $A:=\ms{N}_M(M)$. 
As in the first paragraph of the proof of Theorem~\ref{T:Nashsoft}, we write $V$ as the vanishing set of a single Nash function, 
$V=\{x\in M:\ f(x) = 0\}$ for some $f\in A$. 
Then we denote the complement of $V$ in $M$ by $E(f)$. 
By using Lemma~\ref{L:homeomorphism}, we can identify $E(f)$ with the intersection of the space of maximal ideals of $A$ with a basic Zariski-open subset,
$$
\beta(E(f))= \beta(M)\cap D(f),\quad\text{ where }\quad D(f) :=\{ y\in \text{Nsp}(M):\ f(y)\neq 0\}.
$$  
In particular, $\{E(f):\ f\in A\}$ gives a basis for the $\mathfrak{N}$-topology on $M$. 
Hence, every sheaf on $M$ is characterized by its values on the sets $E(f)$, $f\in A$. 
\begin{definition}
A sheaf of $\ms{N}_M$-modules $\ms{G}$ is called a quasi-coherent sheaf of $\ms{N}_M$-modules if there exists an $A$-module $K$ such that for every $f\in A$ we have 
$$
\ms{G}(E(f))= A_f:=\{ g/f^m :\ \text{$g\in A$ and $m\in \N$}\}.
$$
\end{definition}

\begin{theorem}\label{t:restriction}
We maintain the notation from the above paragraphs. 
In particular, $\beta: M\to \text{Nsp}(M)$ denotes the homeomorphism that is defined in Lemma~\ref{L:homeomorphism}. 
Then the direct image functor $\beta_*$ gives an equivalence between the category of quasi-coherent $\ms{N}_M$-modules and the category of quasi-coherent $\text{Nsp}(M)$-modules.
\end{theorem}

\begin{proof}
Since $\beta(M)$ is precisely the (dense) subset of closed points of $\text{Nsp}(M)$, the Zariski open subsets of $\text{Nsp}(M)$ are uniquely determined by the open Nash subsets of $M$. 
For every pair of open subsets $U\subset V$ in $\text{Nsp}(M)$, we have a commuting diagram of maps as in Figure~\ref{F:comdiagram},
where the vertical maps are restrictions, and the vertical maps are given by $\beta$.
\begin{figure}[htp]
\begin{center}
\begin{tikzcd}
\ms{N}_M(\beta^{-1}(V)) \arrow[r] \arrow[d]
&  \ms{O}_{\text{Nsp}(M)}(V) \arrow[d] \\
\ms{N}_M(\beta^{-1}(U)) \arrow[r]
& \ms{O}_{\text{Nsp}(M)}(U),
\end{tikzcd}
\end{center}
\caption{}
\label{F:comdiagram}
\end{figure}
Hence, $\beta$ defines a map of locally ringed spaces. 
Let $x\in M$. 
The stalk of $\ms{N}_M$ at $x$ and the stalk of $\ms{O}_{\text{Nsp}(M)}$ at $\beta(x)$ are the same. 
Indeed, the stalk $\ms{N}_{M,x}$ is a local ring and the ideal $I_x$ consisting of Nash functions vanishing at $x$ is a maximal ideal by Lemma~\ref{L:closedpoints}.
Therefore, the stalk $\ms{N}_{M,x}$ is the localization of $A:=\ms{N}_M(M)$ at $I_x$. 
But this is precisely the stalk of the structure sheaf $\ms{O}_{\text{Nsp}(M)}$ of the affine scheme $\text{Spec} (A)$ at the point $I_x$.
After this observation, the fact that $(\beta,\beta^\#) : (M,\ms{N}_M)\to (\beta(M),\ms{O}_{\text{Nsp}(M)} \vert_{\beta(M)})$ are isomorphic as locally ringed spaces is clear. Thus, $\beta_*$ defines a functor between the category of sheaves on $M$ and the category of sheaves on $\text{Nsp}(M)$.

Let $K$ be an $A$-module. 
Then $K$ defines a quasi-coherent sheaf of $\ms{O}_{\text{Nsp}(M)}$-modules $\ms{F}$ as well as a quasi-coherent sheaf of $\ms{N}_M$-modules $\ms{G}$. 
The stalk of the direct image $\beta_*\ms{G}$ at $\beta(x)$, where $x\in M$, is isomorphic to the stalk of the restriction of $\ms{F}$ at $\beta(x)$.
This is a consequence of the following facts: 1) there is a 1-1 correspondence between the open subsets of $\beta(M)$ and $\text{Nsp}(M)$, 2) the spaces of global sections of $\ms{F}$ and $\ms{G}$ on the basic open subsets $D(f)$ and $E(f)$, where $f\in A$, are the same $A$-module, $A_f$. 
Hence, we see that $\beta_*\ms{G} = \ms{F}$. 
This means that $\beta_*$ restricts to give a functor from the category of quasi-coherent sheaves of $\ms{N}_M$-modules to the category of quasi-coherent sheaves of $\ms{O}_{\text{Nsp}(M)}$-modules. By abuse of notation, we will denote this subfunctor by $\beta_*$, also. 
We proceed to show that this restricted functor is an equivalence of categories. 
In other words, we will show that it is a full, faithful, and essentially surjective functor. 

The essentially surjectivity of $\beta_*$ readily follows from our arguments in the first paragraph of this proof. 
To prove that $\beta_*$ is full and faithful, we will show that it maps bijectively the short exact sequences of quasi-coherent $\ms{N}_M$-modules to the short exact sequences of quasi-coherent $\ms{O}_{\text{Nsp}(M)}$-modules. 
It follows from the first paragraph of this proof that every exact sequence of sheaves of $\ms{N}_M$-modules of the form 
\begin{align}\label{a:initiallyexact}
0\to \ms{F}'\to \ms{F} \to \ms{F}''\to 0
\end{align} 
gives rise to a short exact sequence of $\ms{O}_{\text{Nsp}(M)} \vert_{\beta(M)}$-modules of the form
\begin{align}\label{a:initiallyexact2}
0\to \beta_*\ms{F}'|_{\beta(M)} \to \beta_*\ms{F}|_{\beta(M)} \to \beta_*\ms{F}''|_{\beta(M)}\to 0.
\end{align}
Assuming that (\ref{a:initiallyexact2}) is a short exact sequence of quasi-coherent $\ms{O}_{\text{Nsp}(M)}|_{\beta(M)}$-modules,
we proceed to show that 
\begin{align}\label{a:quasicoherentNsp}
0\to \beta_*\ms{F}' \to \beta_*\ms{F} \to \beta_*\ms{F}''\to 0
\end{align}
is a short exact sequence of quasi-coherent $\ms{O}_{\text{Nsp}(M)}$-modules. 
Equivalently, we will show that for every point $P\in \text{Nsp}(M)$ the corresponding sequence of stalks, 
\begin{align*}
0\to  (\beta_*\ms{F}')_P \to (\beta_*\ms{F})_P \to (\beta_*\ms{F}'')_P\to 0,
\end{align*}
is exact. 
Recall that two modules are isomorphic if and only if their localizations at every maximal ideal are isomorphic.
Hence, we have to check that 
\begin{align}\label{a:restrictionisexact}
0\to  (\beta_*\ms{F}')_I \to (\beta_*\ms{F})_I \to (\beta_*\ms{F}'')_I\to 0
\end{align}
is exact for every maximal ideal $I\subseteq A$.
But this was our initial assumption in (\ref{a:initiallyexact2}).
Hence, (\ref{a:quasicoherentNsp}) is a short exact sequence of quasi-coherent $\ms{O}_{\text{Nsp}(M)}$-modules. 
Clearly, these implications are reversible. 
In other words, if we start with a short exact sequence of quasi-coherent sheaves of $\ms{O}_{\text{Nsp}_M}$-modules, then its restriction to $\beta(M)$ gives a short exact sequence of quasi-coherent $\beta_*\ms{N}_M$-modules. 
Since $(\beta,\beta^\sharp)$ is an isomorphism, we obtain a short exact sequence of quasi-coherent $\ms{N}_M$-modules. 
Hence, the proof of our theorem is complete. 
\end{proof}

We are now ready to finish the proof our Theorem~\ref{T3:intro} from the introduction.
We record it here as a corollary of Theorem~\ref{t:restriction}.
\begin{corollary}(Theorem~\ref{T3:intro})
Let $M$ be an affine Nash manifold. 
If $\ms{G}$ is a quasi-coherent sheaf of $\ms{N}_M$-modules, then we have $H^i(M,\ms{F}) = 0$ for every $i>0$.  
\end{corollary}
\begin{proof}
Let $\mathfrak{Qcn}(M)$ (resp. $\mathfrak{Qco}(\text{Nsp}(M))$) denote the category of quasi-coherent sheaves of $\ms{N}_M$-modules 
(resp. the category of quasi-coherent sheaves of $\ms{O}_{\text{Nsp}(M)}$-modules). 
By Theorem~\ref{t:restriction}, $\ms{G}$ gives a quasi-coherent sheaf of $\ms{O}_{\text{Nsp}(M)}$-modules, $\ms{F}:=\beta_*\ms{G}$, 
where $\beta_*:\mathfrak{Qcn}(M)\to \mathfrak{Qco}(\text{Nsp}(M))$ is the equivalence of categories.  
We know from~\cite[Ch III, Theorem 3.5]{Hartshorne} that 
$$H^i(\text{Nsp}(M),\ms{F})=0\quad \text{ for every $i>0$}.$$ 
Since $\mathfrak{Qco}(\text{Nsp}(M))$ is equivalent to $\mathfrak{Qcn}(M)$, the vanishing of the higher cohomology of $\ms{F}$ implies the vanishing of the higher cohomology of $\ms{G}$. This finishes the proof of our assertion. 
\end{proof}

Before revisiting the definition of an affine Nash vector bundle, let us briefly discuss Nash groups. In the preliminaries section, we introduced affine semialgebraic groups as the group objects within the category of affine semialgebraic spaces. Similarly, we establish a parallel definition for Nash groups.

\begin{definition}
A {\em Nash group} is a group which is at the same time a Nash manifold such that all group operations are Nash maps. 
A {\em Nash homomorphism} is a group homomorphism between two Nash groups which is at the same time a Nash map. 
If the underlying Nash manifold of a Nash group $G$ is an affine Nash manifold, then 
we call $G$ an {\em affine Nash group}. 
\end{definition}
In other words, (affine) Nash groups are precisely the group objects in the category of 
(affine) Nash manifolds. A general introduction to Nash groups can be found in Sun's article~\cite{Sun2015}.
\medskip

Let $K$ be a real algebraic group. 
Then the underlying real algebraic variety of $K$ is a smooth manifold.
Furthermore, the group operations of $K$ are smooth maps. 
It follows that every real algebraic group is naturally a Nash manifold. 
Here is an example of an affine Nash group which is not a real algebraic group.

\begin{example}\label{E:firstH}
Let $H$ denote the group of positive real numbers under multiplication. 
Then $H$ is an affine Nash group. 
We will show that $H$ is not a real algebraic group. 
To this end, we consider the $\R$-algebra of regular functions on $H$, denoted by $\R[H]$. 
It is clear that $\R[H]$ contains all rational functions of the form $\frac{1}{x+a}$, where $a\in H$. 
These rational functions are all algebraically independent over $\R[x]$, which implies that $\R[H]$ is not a finitely generated $\R$-algebra. 
Therefore, $H$ is not a real algebraic variety. 
In particular, it follows that $H$ is not a real algebraic group. 
\end{example}

We will use our previous example to illustrate several subtle distinctions between the theory of real algebraic groups and the theory of complex algebraic groups.

First, we will show that the connected component of the identity element in a real algebraic group need not be a real algebraic group. This is in contrast to the case of complex algebraic groups.
\medskip

Let $G$ denote the set $G:=\{(x,y) \in \R^2 :\ xy = 1 \}$, which is a real algebraic group with respect to entry-wise multiplication. 
Geometrically, $G$ is a hyperbola in the real plane. 
The defining ideal of $G$ is the prime ideal $(xy-1)$ in $\R[x,y]$.
Therefore, $G$ is irreducible as a real algebraic variety. 
However, it has two connected components in the ordinary topology. 
The connected component containing the identity element $(1,1)\in G$ is Nash isomorphic to $H$ via the Nash imbedding 
\begin{align*}
f: H &\longrightarrow G,\\ 
x &\longmapsto (x,x^{-1}).
\end{align*}
We conclude from this example that the connected component of the identity element in a real algebraic group need not be a real algebraic subgroup. 

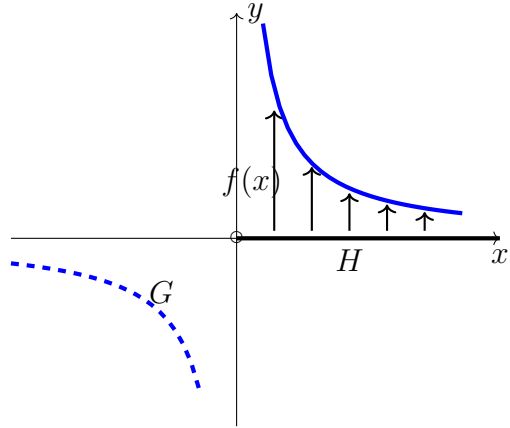
\begin{center}
\begin{tikzpicture}
\draw[->] (-3,0) -- (3.5,0) node[below] {$x$};
\draw[->] (0,-2.5) -- (0,3) node[right] {$y$};
\draw[ultra thick, -] (0,0) -- (3.5,0);

\draw[thick, ->] (0.5,.1) -- (0.5,1.7);
\draw[thick, ->] (1,.1) -- (1, .95);
\draw[thick, ->] (1.5,.1) -- (1.5, .6);
\draw[thick, ->] (2,.1) -- (2, .45);
\draw[thick, ->] (2.5,.1) -- (2.5, .35);

\draw[blue,ultra thick, domain=0.35 : 3] plot (\x, 1/\x);
\draw[blue,ultra thick, dashed, domain=-3:- 0.5] plot (\x, 1/\x);

\node at (0,0) {$\circ$};
\node[below] at (1.5,0) {$H$};
\node[above] at (-1,-1) {$G$};
\node[left] at (.75,.75) {$f(x)$};

\end{tikzpicture}
\end{center}

Next, we will show that even a homomorphic image of a real algebraic group under a real algebraic group homomorphism need not be a real algebraic group. This is in contrast to the case of complex algebraic groups.
Once again, we focus on our previous example. 

\begin{example}
Let $G$ denote the multiplicative group of nonzero real numbers, $(\R^*,\times)$.
Then $G$ is isomorphic to the hyperbola of the previous example. 
In particular, $G$ is a real algebraic group.
Since $G$ is abelian, the squaring map, $sq: G\to G$, $g\mapsto g^2$ is an algebraic group homomorphism. 
The image of $sq$ is canonically isomorphic to the group $(\R^+,\times)$, which we denoted by $H$ in Example~\ref{E:firstH}. 
As we showed in that example, $H$ is not a real algebraic group but an affine Nash group.
\end{example}

An important result of Hrushovsky and Pillay~\cite{HrushovskiPillay1994} states 
that for every connected affine Nash group $H$, there exists a complex algebraic group $G$, defined over $\R$, that admits a finite surjective Nash group homomorphism $H\to G(\R)^0$, where $G(\R)^0$ is the connected component of the $\R$-rational points of $G$ such that $1\in G(\R)^0$.
In other words, affine Nash groups are precisely the finite covers of real algebraic groups. 
\medskip

We go back to our discussion about the $\mathfrak{N}$-topology.
We will consider the affine Nash vector bundles in our topology.  
A locally free sheaf of $\ms{N}_M$-modules $\ms{E}$ on $M$ is said to be an affine Nash vector bundle of class $\ms{N}_M$ 
if there is a finite open cover $M=\bigcup_{i=1}^m U_i$ in the $\mathfrak{N}$-topology such that 
for every $i\in \{1,\dots, m\}$ there is an isomorphism $\tau_i : \mc{E} |_{U_i} \to \ms{N}_M |_{U_i} \otimes_\R D$, 
and the transition functions $g_{i,j}:= \tau_i \circ \tau_j^{-1}$ ($i,j\in \{1,\dots, m\}$) are isomorphisms of $\ms{N}_M$-modules. 
We note that the maps induced by the transition functions on the stalks $\ms{N}_{M,x}\otimes D$ are contained in $\text{GL}(D)$. 
The structure group of an affine Nash vector bundle is naturally an affine Nash group.
\medskip

Let us recall the Nash analog of a ringed $D$-space structure.
A sheaf of $k$-algebras $\ms{A}$ on an affine Nash manifold $M$ is called a {\em ringed $D$-space structure on $M$} if there exists a map of sheaves of $k$-algebras 
$\vep : \ms{A} \to \ms{N}_M$, and a covering $\{ U_i \to X | i\in I \}$ in the $\mathfrak{N}$-topology equipped with the maps of sheaves of $\ms{N}_M$-algebras, 
\begin{align*}
\tau_i : \ms{A} \vert_{U_i} \longrightarrow (\ms{N}_M\otimes D) \vert_{U_i}\qquad (i\in I) 
\end{align*}
such that $(id \otimes \vep_0) \circ \tau_i = \vep$.
The hypothesis of our next proposition is a specialization of the hypothesis of Theorem~\ref{T:maintool}.
Hence we have the corresponding conclusion.
We omit the details of its proof.

\begin{lemma}\label{L:BR_Nash}
Let $(M,\ms{N}_M)$ be an affine Nash manifold.
Let $\vep_\circ: D\to \R$ be a surjective homomorphism of finite dimensional $\R$-algebras such that $J:=\ker \vep_\circ$
is a nilpotent ideal. 
Let $\vep : \ms{A}\to \ms{N}_M$ be an affine Nash $D$-manifold structure on $M$.
Let $\ms{I}$ denote the kernel of $\vep$. 
If the sheaf defined by $\ms{I}_0:=\ms{I} \cap Z(\ms{A})$ satisfies the conditions 
\begin{align*}
H^1 (M, \text{Der}(\ms{N}_M, \ms{I}_0^p / \ms{I}^{p+1}_0))=0
\end{align*} 
for every $p\in \Z_+$, then there is a map of sheaves of $\ms{N}_M$-algebras $\varphi: \ms{N}_M\to Z(\ms{A})$
splitting $\vep$. 
In this case, as a sheaf of $\ms{N}_M$-modules, $\ms{A}$ is an affine Nash vector bundle over $M$ of class $\ms{N}_M$ with structure group in $\text{Aut}(D,\vep_0)$. 
\end{lemma}

For the rest of this section, we will denote by $D$ an exterior algebra $\LL \R^s$ ($s\in \Z_+$) with grading as described at the beginning of  Section~\ref{S:Locally}. 
Let us recall the statement of Theorem~\ref{T2:intro} from the introduction section.
\medskip

Let $(M,\ms{A})$ be an affine Nash supermanifold of odd-dimension $s$. 
Then there exists a Nash vector bundle $E$ on $M$ and a finite cover $M=\bigcup_{i\in I} U_i$ by $\mathfrak{N}$-open sets such that, 
for each $i\in I$, there is an isomorphism of $\Z_2$-graded algebras, $\ms{A}(U_i) \cong \Gamma( U_i, \LL E)$.

\begin{proof}[Proof of Theorem~\ref{T2:intro}]
We begin with showing the vanishing of all higher cohomologies of the sheaves $\text{Der}(\ms{N}_M, \ms{I}_0^p / \ms{I}^{p+1}_0)$ for $p\in \Z_+$.
To this end, we notice that $\text{Der}(\ms{N}_M, \ms{I}_0^p / \ms{I}^{p+1}_0)$ is a quasi-coherent sheaf of $\ms{N}_M$-modules. 
Hence, by Theorem~\ref{T3:intro}, we have 
\begin{align*}
H^i ( M,\text{Der}(\ms{N}_M, \ms{I}_0^p / \ms{I}^{p+1}_0))=0 \quad \text{ for every $\{i,p\}\subset \Z_+$}.
\end{align*} 
Now, our Lemma~\ref{L:BR_Nash} implies that, for every $D$-manifold structure $\vep: \ms{A}\to \ms{N}_M$, there is a splitting homomorphism $\varphi: \ms{N}_M\to Z(\ms{A})$. 
Hence, we can view our sheaf of $\Z_2$-graded commutative algebras $\ms{A}$ on $M$ as a vector bundle over $M$ of class $\ms{N}_M$
with structure group in $\text{Aut}_0(D,\vep_0)$. 
Note here that $\text{Aut}_0(D,\vep_0)$ is the subgroup of $\text{Aut}(D,\vep_0)$ consisting of $\Z_2$-grading preserving elements. 
Since both of these groups are real algebraic groups, we can view them as affine Nash groups. 
To finish the proof, we repeat the reasoning of the second half of the proof of our Theorem~\ref{T1:intro}. 
Since the arguments are almost identical except that we now have affine Nash subgroups instead of affine semialgebraic groups, 
we conclude that the structure group of $\ms{A}$ can be reduced to $\text{GL}_s(\R)$ in the category of affine Nash groups. 
This finishes the proof of our second main result. 
\end{proof}

\section{Final Remarks}\label{S:Final}

Every Nash group $G$ is automatically a real Lie group since the underlying space of $G$ is a smooth manifold and the group operations are $\ms{C}^\infty$ maps. As we mentioned earlier, a theorem of Malgrange shows that these maps are analytic maps. 
The classification of one-dimensional (real) Lie groups is a straightforward process. 
Indeed, a connected one-dimensional Lie group is isomorphic to either $(\R,+)$ or $(\mathbb{S},\cdot)$, where 
$\mathbb{S}:=\{ z\in \C^* :\ |z|=1\}$, 
and the group operation of $\mathbb{S}$ is the multiplication of complex numbers. 
However, the classification of connected one-dimensional Nash groups turned out to be substantially more intricate.
This task was successfully accomplished by Madden and Stanton in their work~\cite{MaddenStanton1992}.
A readily accessible exposition of the results of Madden and Stanton can be found in the comprehensive article~\cite{Shiota1996}.

For the classification of connected one-dimensional Nash groups, an essential ingredient is the notion of a {\em locally Nash group}.
Such a structure can be defined as a group object in the category of {\em locally Nash manifolds}.
Here, a locally Nash manifold is defined as an analytic manifold with a possibly infinite system of Nash coordinate neighborhoods. 
Equivalently, a {\em locally Nash manifold} $M$ is a group object in the category of locally semialgebraic spaces such that the underlying space of $M$ is an analytic manifold, and the group operations of $M$ are analytic maps.  

\begin{example}
Let us consider the unit circle $\mathbb{S}$ (with its group structure inherited from $\C^*$) as an affine Nash group.
Let $\widetilde{\mathbb{S}}$ denote the universal covering space of $\mathbb{S}$. 
Then $\widetilde{\mathbb{S}}$ has the structure of a locally Nash manifold and a compatible group structure on it. 
In other words, $\widetilde{\mathbb{S}}$ is locally Nash group. 
Then $\widetilde{\mathbb{S}}$ is not a Nash group. 
To see this, consider the covering homomorphism $p: \widetilde{\mathbb{S}} \to \mathbb{S}$.
If $\widetilde{\mathbb{S}}$ was a Nash group, then $p$ would be a Nash group homomorphism. 
However, the kernel of $p$ is isomorphic to $\Z$, which is not a Nash group but a locally Nash group. 
It follows that $\widetilde{\mathbb{S}}$ is not a Nash group. 
\end{example}

We will now provide a very brief overview of the development of the classification of connected one-dimensional Nash groups.
\medskip

To begin with, it can be shown that the universal covering group of a locally Nash group is also a locally Nash group. This observation leads to a separation of the classification problem into two distinct parts:
\begin{enumerate}
\item The classification of isomorphism classes of simply connected locally Nash groups.
\item The classification of isomorphism classes of the quotients of a simply connected locally Nash group. 
\end{enumerate}
Let $G$ be a one-dimensional, connected, and simply connected locally Nash group.
Then the underlying locally Nash manifold is locally Nash isomorphic to $\R$ with an appropriate (possibly infinite) system of locally Nash coordinate charts. 
In particular, there is an analytic Nash coordinate chart $\psi : U\to \R$ such that $0\in U$, and the Nash group operations on every sufficiently small neighborhood of 0 are transferred via $\psi$ from $(\R,+)$. 
Furthermore, this coordinate chart $\psi$ determines, up to a locally Nash group isomorphism, the locally Nash group structure on the whole of $G$. Therefore, it makes sense to denote the locally Nash group $G$ by the pair $(\R,\psi)$. 
The gist of the classification of the pairs $(\R,\psi)$ is the {\em addition theorem of Weierstrass}. All of this is explained in more detail in the original article~\cite{MaddenStanton1992} as well as in the Shiota's article~\cite[pages 107--108]{Shiota1996}. 
\medskip

In the present article, we have introduced the categories of Nash supermanifolds and locally semialgebraic superspaces implicitly.
Therefore, the following definition now makes sense.  
\begin{definition}
A {\em Nash supergroup} (or a {\em $\Z_2$-graded Nash group}) is a group object in the category of Nash supermanifolds.
Likewise, a {\em locally semialgebraic supergroup} (or a {\em $\Z_2$-graded locally semialgebraic group}) is a group object in the category of locally semialgebraic supergroups. 
\end{definition} 
Continuing our work in this paper and building upon the research of Madden and Stanton, we have classified the connected one-dimensional affine locally Nash supergroups. The details of this work will appear elsewhere.

\section*{Acknowledgements}
Final edits of this paper were completed while we were visiting the Okinawa Institute of Science and Technology (OIST) through the Theoretical Sciences Visiting Program (TSVP). We thank OIST-TSVP for providing us with a comfortable work environment. 
We gratefully acknowledge the support of the Louisiana Board of Regents (grant no. 090ENH-21) which partially funded this work. 
Additionally, we extend our gratitude to the anonymous referee for their meticulous review of our manuscript, offering invaluable feedback that greatly enhanced the quality of our article. Our thanks also go to J\"org Feldvoss and Naufil Sakran for their insightful review, which significantly contributed to the clarity of our writing.

\bibliography{references}

\providecommand{\bysame}{\leavevmode\hbox to3em{\hrulefill}\thinspace}
\providecommand{\MR}{\relax\ifhmode\unskip\space\fi MR }
\providecommand{\MRhref}[2]{%
  \href{http://www.ams.org/mathscinet-getitem?mr=#1}{#2}
}
\providecommand{\href}[2]{#2}
\begin{thebibliography}{{Sta}22}

\bibitem[AG08]{AizenbudGourevitch2008}
Avraham Aizenbud and Dmitry Gourevitch, \emph{Schwartz functions on {N}ash
  manifolds}, Int. Math. Res. Not. IMRN (2008), no.~5, Art. ID rnm 155, 37.
  \MR{2418286}

\bibitem[Art62]{Artin}
Michael Artin, \emph{Grothendieck topologies: notes on a seminar}, Harvard
  University, Dept. of Mathematics, 1962.

\bibitem[Bat79]{Batchelor1979}
Marjorie Batchelor, \emph{The structure of supermanifolds}, Trans. Amer. Math.
  Soc. \textbf{253} (1979), 329--338. \MR{536951}

\bibitem[Bat80]{Batchelor1980}
\bysame, \emph{Two approaches to supermanifolds}, Trans. Amer. Math. Soc.
  \textbf{258} (1980), no.~1, 257--270. \MR{554332}

\bibitem[BCR98]{BochnakCosteRoy}
Jacek Bochnak, Michel Coste, and Marie-Fran\c{c}oise Roy, \emph{Real algebraic
  geometry}, Ergebnisse der Mathematik und ihrer Grenzgebiete (3) [Results in
  Mathematics and Related Areas (3)], vol.~36, Springer-Verlag, Berlin, 1998,
  Translated from the 1987 French original, Revised by the authors.
  \MR{1659509}

\bibitem[BL75]{BerezinLeites}
F.~A. Berezin and D.~A. Le\u{\i}tes, \emph{Supermanifolds}, Dokl. Akad. Nauk
  SSSR \textbf{224} (1975), no.~3, 505--508. \MR{0402795}

\bibitem[BR84]{BlattnerRawnsley}
Robert~J. Blattner and John~H. Rawnsley, \emph{Remarks on {B}atchelor's
  theorem}, Mathematical aspects of superspace ({H}amburg, 1983), NATO Adv.
  Sci. Inst. Ser. C: Math. Phys. Sci., vol. 132, Reidel, Dordrecht, 1984,
  pp.~161--171. \MR{773083}

\bibitem[Del91]{Delfs}
Hans Delfs, \emph{Homology of locally semialgebraic spaces}, Lecture Notes in
  Mathematics, vol. 1484, Springer-Verlag, Berlin, 1991. \MR{1176311}

\bibitem[DeW84]{DeWitt}
Bryce DeWitt, \emph{Supermanifolds}, Cambridge Monographs on Mathematical
  Physics, Cambridge University Press, Cambridge, 1984. \MR{778559}

\bibitem[DK81]{DelfsKnebusch1981II}
Hans Delfs and Manfred Knebusch, \emph{Semialgebraic topology over a real
  closed field. {II}. {B}asic theory of semialgebraic spaces}, Math. Z.
  \textbf{178} (1981), no.~2, 175--213. \MR{631628}

\bibitem[DK84]{DelfsKnebusch1984}
\bysame, \emph{An introduction to locally semialgebraic spaces}, vol.~14, 1984,
  Ordered fields and real algebraic geometry (Boulder, Colo., 1983),
  pp.~945--963. \MR{773141}

\bibitem[DK85]{DelfsKnebusch}
\bysame, \emph{Locally semialgebraic spaces}, Lecture Notes in Mathematics,
  vol. 1173, Springer-Verlag, Berlin, 1985. \MR{819737}

\bibitem[DM99]{DeligneMorgan1999}
Pierre Deligne and John~W. Morgan, \emph{Notes on supersymmetry (following
  {J}oseph {B}ernstein)}, Quantum fields and strings: a course for
  mathematicians, {V}ol. 1, 2 ({P}rinceton, {NJ}, 1996/1997), Amer. Math. Soc.,
  Providence, RI, 1999, pp.~41--97. \MR{1701597}

\bibitem[Efr74]{Efroymson1974}
Gustave~A. Efroymson, \emph{A {N}ullstellensatz for {N}ash rings}, Pacific J.
  Math. \textbf{54} (1974), 101--112. \MR{360576}

\bibitem[Efr82]{Efroymson}
\bysame, \emph{The extension theorem for {N}ash functions}, Real algebraic
  geometry and quadratic forms ({R}ennes, 1981), Lecture Notes in Math., vol.
  959, Springer, Berlin-New York, 1982, pp.~343--357. \MR{683141}

\bibitem[EP16]{EdmundoPrelli}
M\'{a}rio~J. Edmundo and Luca Prelli, \emph{Sheaves on
  {$\mathcal{T}$}-topologies}, J. Math. Soc. Japan \textbf{68} (2016), no.~1,
  347--381. \MR{3454562}

\bibitem[God58]{Godement}
Roger Godement, \emph{Topologie alg\'{e}brique et th\'{e}orie des faisceaux},
  Actualit\'{e}s Scientifiques et Industrielles [Current Scientific and
  Industrial Topics], No. 1252, Hermann, Paris, 1958, Publ. Math. Univ.
  Strasbourg. No. 13. \MR{0102797}

\bibitem[Har77]{Hartshorne}
Robin Hartshorne, \emph{Algebraic geometry}, Springer-Verlag, New
  York-Heidelberg, 1977, Graduate Texts in Mathematics, No. 52. \MR{0463157}

\bibitem[HP94]{HrushovskiPillay1994}
Ehud Hrushovski and Anand Pillay, \emph{Groups definable in local fields and
  pseudo-finite fields}, Israel J. Math. \textbf{85} (1994), no.~1-3, 203--262.
  \MR{1264346}

\bibitem[Knu71]{Knutson}
Donald Knutson, \emph{Algebraic spaces}, Lecture Notes in Mathematics, Vol.
  203, Springer-Verlag, Berlin-New York, 1971. \MR{0302647}

\bibitem[Kos77]{Kostant1975}
Bertram Kostant, \emph{Graded manifolds, graded {L}ie theory, and
  prequantization}, Differential geometrical methods in mathematical physics
  ({P}roc. {S}ympos., {U}niv. {B}onn, {B}onn, 1975), 1977, pp.~177--306.
  Lecture Notes in Math., Vol. 570. \MR{0580292}

\bibitem[Mal67]{Malgrange}
B.~Malgrange, \emph{Ideals of differentiable functions}, Tata Institute of
  Fundamental Research Studies in Mathematics, vol.~3, Tata Institute of
  Fundamental Research, Bombay; Oxford University Press, London, 1967.
  \MR{0212575}

\bibitem[Man97]{Manin}
Yuri~I. Manin, \emph{Gauge field theory and complex geometry}, second ed.,
  Grundlehren der mathematischen Wissenschaften [Fundamental Principles of
  Mathematical Sciences], vol. 289, Springer-Verlag, Berlin, 1997, Translated
  from the 1984 Russian original by N. Koblitz and J. R. King, With an appendix
  by Sergei Merkulov. \MR{1632008}

\bibitem[MS92]{MaddenStanton1992}
James~J. Madden and Charles~M. Stanton, \emph{One-dimensional {N}ash groups},
  Pacific J. Math. \textbf{154} (1992), no.~2, 331--344. \MR{1159515}

\bibitem[Pec85]{Pecker}
Daniel Pecker, \emph{On {E}froymson's extension theorem for {N}ash functions},
  J. Pure Appl. Algebra \textbf{37} (1985), no.~2, 193--203. \MR{796409}

\bibitem[Pil88]{Pillay1988}
Anand Pillay, \emph{On groups and fields definable in {$o$}-minimal
  structures}, J. Pure Appl. Algebra \textbf{53} (1988), no.~3, 239--255.
  \MR{961362}

\bibitem[Ris73]{Risler}
Jean-Jacques Risler, \emph{Sur l'anneau des fonctions de {N}ash globales}, C.
  R. Acad. Sci. Paris S\'{e}r. A-B \textbf{276} (1973), A1513--A1516.
  \MR{318510}

\bibitem[Rob83]{Robson}
Robert Robson, \emph{Embedding semi-algebraic spaces}, Math. Z. \textbf{183}
  (1983), no.~3, 365--370. \MR{706394}

\bibitem[Rog80]{Rogers1980}
Alice Rogers, \emph{A global theory of supermanifolds}, J. Math. Phys.
  \textbf{21} (1980), no.~6, 1352--1365. \MR{574696}

\bibitem[Ser56]{GAGA}
Jean-Pierre Serre, \emph{G\'{e}om\'{e}trie alg\'{e}brique et g\'{e}om\'{e}trie
  analytique}, Ann. Inst. Fourier (Grenoble) \textbf{6} (1955/56), 1--42.
  \MR{82175}

\bibitem[Shi87]{Shiota}
Masahiro Shiota, \emph{Nash manifolds}, Lecture Notes in Mathematics, vol.
  1269, Springer-Verlag, Berlin, 1987. \MR{904479}

\bibitem[Shi96]{Shiota1996}
\bysame, \emph{Nash functions and manifolds}, Lectures in real geometry
  ({M}adrid, 1994), De Gruyter Exp. Math., vol.~23, de Gruyter, Berlin, 1996,
  pp.~69--112. \MR{1440210}

\bibitem[{Sta}22]{stacks-project}
The {Stacks project authors}, \emph{The stacks project},
  \url{https://stacks.math.columbia.edu}, 2022.

\bibitem[Sun15]{Sun2015}
Binyong Sun, \emph{Almost linear {N}ash groups}, Chin. Ann. Math. Ser. B
  \textbf{36} (2015), no.~3, 355--400. \MR{3341164}

\end{thebibliography}

\bibliographystyle{amsalpha}

\end{document}